\documentclass[a4paper]{amsart}

\usepackage{amsmath,amsthm,amsfonts,amssymb,amscd}
\input{xy}
\xyoption{all}

\newenvironment{demo}{\begin{proof}[Proof]}{\end{proof}}
\newtheorem{teo}{Theorem}[section]
\newtheorem{pro}[teo]{Proposition}
\newtheorem*{pro1}{Proposition 4.3'}
\newtheorem{lem}[teo]{Lemma}
\newtheorem{cor}[teo]{Corollary}
\newtheorem{defi}[teo]{Definition}
\newtheorem{rem}[teo]{Remark}
\newtheorem*{exa}{Example}

\newenvironment{sis}{\left\{\begin{aligned}}{\end{aligned}\right.}

\numberwithin{equation}{section}

\renewcommand{\P}{\mathbb{P}}
\renewcommand{\O}{\mathcal{O}}
\newcommand{\C}{\mathcal{C}}
\newcommand{\F}{\mathcal{F}}
\newcommand{\M}{\mathcal{M}}
\renewcommand{\L}{\mathcal{L}}
\newcommand{\E}{\mathcal{E}}
\newcommand{\GG}{\mathcal{G}}
\newcommand{\T}{\mathcal{T}_{a,b}}
\newcommand{\Z}{\mathbb{Z}}

\newcommand{\G}{\mathbb{G}_m}

\newcommand{\Pic}{{\rm Pic}}

\newcommand{\Cl}{{\rm Cl}}

\renewcommand{\S}{{\rm Sym}^{2g+2}(\P^1)}

\newcommand{\B}{\mathbb{B}(2,2g+2)}

\newcommand{\Bsm}{\mathbb{B}_{sm}(2,2g+2)}

\newcommand{\A}{\mathbb{A}(2,2g+2)}

\newcommand{\Asm}{\mathbb{A}_{sm}(2,2g+2)}

\newcommand{\AO}{\mathbb{A}(2,2g+2)\backslash\{0\}}

\newcommand{\BS}{\mathbb{B}(2,6)}

\newcommand{\BSsm}{\mathbb{B}_{sm}(2,6)}

\renewcommand{\H}{\mathcal{H}_g}

\newcommand{\D}{\mathcal{D}_{2g+2}}

\renewcommand{\det}{{\rm det}}

\newcommand{\Gdiv}{\mathcal{G}_{\rm div}}

\begin{document}

\title{Families of hyperelliptic curves}
\author{Sergey Gorchinskiy and Filippo Viviani}
\address{Steklov Mathematical Institute, Gubkina str. 8, 119991 Moscow}
\email{serge.gorchinsky@rambler.ru}
\address{Universita' degli studi di Roma Tor Vergata, Dipartimento di
matematica, via della ricerca scientifica 1, 00133 Rome}
\email{viviani@axp.mat.uniroma2.it}
\thanks{The first author was partially supported by RFFI grants 04-01-00613 and
05-01-00455.}
\maketitle

\section{Introduction}

Throughout this work we deal with a natural number $g\ge 2$ and with an
algebraically closed field $k$ whose characteristic differs from 2.
A hyperelliptic curve of genus $g$ over $k$ is a smooth curve of genus $g$, that
is a double cover of the
projective line $\P^1$. The Riemann-Hurwitz formula implies that this covering
should be ramified at $2g+2$ points. 

Because of this explicit description, hyperelliptic curves have been studied for
a long time from different points of view. Among recent advances, we want to
mention the determination of all the possible automorphism groups of
hyperelliptic curves (see \cite{BS}, \cite{BGG}, \cite{sha}) as well as the
extensive use of the Jacobian of hyperelliptic curves in cryptography (see
\cite{Sch}, \cite{Can}, \cite{Kob}, \cite{Fre},
\cite{Gau}, \cite{Ked}, \cite{Lan}, and the survey paper \cite{JMS}).

In this paper we are interested in the moduli space $H_g$
of hyperelliptic curves and in the moduli stack $\mathcal{H}_g$ of hyperelliptic
curves, whose definitions we are going to briefly recall now.

The MODULI SCHEME $H_g$ of hyperelliptic curves is defined as
\begin{equation*}
H_g=(\S-\Delta)/PGL_2,
\end{equation*}
where $\S$ is the $(2g+2)$-th symmetric power of $\P^1$, $\Delta$ is
the closed subset where at the least two points coincide and the
action of $PGL_2$ comes from its natural action on $\P^1$. 
Since a hyperelliptic curve over $k$ is completely determined (up to
isomorphism) by $2g+2$ points on $\P^1$ (up
to isomorphism), over which the corresponding double cover of $\P^1$ ramifies,
$H_g$ has the property
that its closed points parameterize isomorphism classes of hyperelliptic curves.

This modular variety has been studied from different points of 
view:  Katsylo and Bogomolov proved
its rationality (see \cite{Kat}, \cite{Bog}), Avritze and Lange considered
various
compactifications of $H_g$ (that is an affine variety) comparing them with each
other (see \cite{AL}).

Our new contribution to the study of $H_g$ is the determination of the Picard 
group 
$\Pic(H_g)$ and of the divisor class group $\Cl(H_g)$. 
We prove that, away from some bad characteristic of the base 
field, $\Pic(H_g)$ is trivial (theorem \ref{picmoduli2}) 
while $\Cl(H_g)$ is a cyclic group of order $4g+2$ if 
$g\geq 3$ and $5$ if $g=2$ (theorem \ref{picmoduli}).
The fact that $\Pic(H_g)\neq \Cl(H_g)$ indicates that $H_g$ is a singular 
variety 
(although its explicit description as quotient imply that it's a normal variety)  
and in fact we determine its smooth locus in proposition \ref{smooth}. 

The MODULI FUNCTOR $\mathcal{H}_g$ of hyperelliptic curves is the contravariant
functor
$$\mathcal{H}_g : Sch_{/k}\rightarrow Set$$
which associates to every $k$-scheme $S$ the set
$$\mathcal{H}_g(S)=\left\{
\F\rightarrow S \text{ family of hyperelliptic smooth curves of genus } g
\right\}_{/\cong}.$$

Arsie and Vistoli (see \cite{AV} and also \cite{Vis} only for $g=2$) 
proved that $\H$ is a Deligne-Mumford stack isomorphic 
to a quotient stack, precisely  
$$\H=[\Asm/(GL_2/\mu_{g+1})]$$
where $\Asm$ is 
the space of binary forms in two variables of degree $2g+2$ 
having only simple factors and $GL_2/\mu_{g+1}$ acts as 
$[A]\cdot f(x)=f(A^{-1}\cdot x)$. Moreover they compute the 
equivariant Picard group $\Pic^{GL_2/\mu_{g+1}}(\Asm)$ which in fact 
is isomorphic to the Picard  group of the stack $\H$ (as defined functorially 
by Mumford in \cite{Mum}). In the case $g=2$, Vistoli proved (in \cite{AV})
that this group is generated by the first Chern class of the Hodge bundle but in   
the general case there wasn't known such a functorial description.

In theorem \ref{gen}, we provide an explicit functorial description of a 
generator 
of the Picard group of the stack. Moreover in theorem $\ref{elements}$
we consider natural elements of the Picard group (obtaining by pushing-forward
linear combinations of the relative canonical divisor and the relative 
Weiertrass
divisor and then taking the determinant) and express them in terms of the 
generator found above. In particular we prove that the first Chern class
of the Hodge bundle generates the Picard group if and only if $4$ doesn't 
divide $g$ in which case it generates a subgroup of index $2$ 
(corollary \ref{Hodge}). 

It is well known that $H_g$ is a COARSE MODULI SCHEME for the functor
$\mathcal{H}_g$, which means that there is a natural transformation of functors
$$\Phi_{\mathcal{H}}:\mathcal{H}_g\rightarrow {\rm Hom}(-,H_g)$$
satisfying two properties:
\begin{itemize}
\item[(i)] $\Phi_{\mathcal{H}}({\rm Spec}(K)):\mathcal{H}_g({\rm 
Spec}(K))\rightarrow
{\rm Hom}({\rm Spec}(K),H_g)$ is bijective for any algebraically closed field
$K=\overline{K}\supset k$.
\item[(ii)] (Universal property) If $N$ is a scheme and
$\Psi:\mathcal{H}_g\rightarrow {\rm Hom}(-,N)$
is a natural transformation of functors, then there exists a unique morphism
$\pi:H_g\rightarrow N$
such that the corresponding natural transformation  $\Pi:{\rm Hom}(-
,H_g)\rightarrow
{\rm Hom}(-,N)$
satisfies $\Psi=\Pi \circ \Phi_{\mathcal{H}}$.
\end{itemize}

Another problem we treat in this work is the question: 
is $H_g$ a FINE moduli scheme for the functor $\mathcal{H}_g$? And if not,
how far is from being such? By definition, being a fine moduli scheme 
would mean that $\Phi_{\mathcal{H}}$ is an isomorphism of functors, 
or, in other words, that
there exists a universal family of hyperelliptic curves
$\mathcal{F}_g\rightarrow H_g$ such that every other family
$f:\mathcal{F}\rightarrow S$ is
obtained from this one by pulling back via the modular map
$\Phi_S(f):S\rightarrow H_g$, i.e.
$\mathcal{F}\cong \mathcal{F}_g \times_{H_g} S$.

To attack this problem, we introduce a new moduli functor $\D$ which 
is intermediate between $\H$ and $H_g$ and is defined as the contravariant
functor
$$\D : Sch_{/k}\rightarrow Set$$
which associates to every $k$-scheme $S$ the set
$$\mathcal{D}_{2g+2}(S)=\left\{
\begin{aligned}
&\C\rightarrow S \text{ family of } \mathbb{P}^1 \text{ and }
D\subset \C
\text{ an effective Cartier divisor }\\
&\text{ finite and \'etale over } S \text{ of degree }
2g+2
\end{aligned}
\right\}_{/\cong}$$

Since, by general results of Lonstead and Kleiman (\cite{LK}), a hyperelliptic 
family 
is a double cover of a family of $\P^1$ ramified along a Cartier divisor $D$ as 
in the definition
above, there is a natural transformation of functors $\Phi:\H\to \D$. Moreover,
since over an algebraically closed field giving a hyperelliptic curve $C$ 
is equivalent to give the $(2g+2)$-points (up to isomorphism) where the $2:1$ 
map  
$C\to \P^1$ ramifies, both these moduli functor have $H_g$ as coarse moduli 
scheme.  

We prove (theorem \ref{D-stack}) that $\D$ is actually an algebraic stack 
isomorphic to a quotient stack, precisely
$$\D=[\Bsm/(PGL_2)],$$
where $\Bsm$ is the projective space of smooth binary forms in $2$ variables of 
degree
$2g+2$ and the action of $PGL_2$ is defined by $[A]\cdot [f(x)]=[f(A^{-1}\cdot 
x)]$.
Using this description as a stack, we compute the Picard group of $\D$ (giving 
an 
explicit generator) and prove that the natural pull-back map 
$\Pic(\D)\to \Pic(\H)$ is an isomorphism for $g$ even and an injection of index 
$2$ 
for $g$ odd (theorem \ref{D-pic}).

Next, after this digression into the study of the auxiliary functor $\D$, we 
return to the 
study of the finess of $H_g$ for $\H$ and $\D$. Since the existence  
of automophisms is always one of the most serious obstructions to the finess of 
a moduli scheme, we restrict to the open subset $H_g^0$ of hyperelliptic 
curves with automorphism group reduced to the hyperelliptic involution
(which we call hyperelliptic curves without extra-automorphisms) as well
as to the corresponding open substack $\H^0:=\H\times_{H_g} H_g^0$ 
and $\D^0=\D\times_{H_g} H_g^0$.   

The first result is that $H_g^0$ is actually a fine moduli scheme for the 
functor $\D^0$, that is over $H_g^0$ there exists a universal family 
of $\P^1$ together with a universal Cartier divisor $D$ as above (theorem 
\ref{universaldivisor}). On the other hand, there doesn't exist over $H_g^0$
a universal family of hyperelliptic curves, thus $H_g^0$ is not fine for 
$\H^0$. In fact, in theorem \ref{doublecovers},  we prove that
the set of families of hyperelliptic curves (without extra-automorphisms) over 
$S$ 
with a fixed modular map (if non empty) is a principal homogeneous space for 
$H^1_{\acute e t}(S,\Z/2\Z)$.

Then we deal with the existence
of a tautological family of hyperelliptic curves over an open subset of $H_g$ 
and
we prove that a tautological family exists over an open subset 
if and only if $g$ is odd (theorem \ref{tautologicalfamily}). 
This was stated in the exercise 2.3 of the book 
of Harris-Morrison (\cite{HM}) but with a mistake: they say universal family 
but in fact it's only a tautological family by what said before!
 
The paper is organized as follows. In section $2$ we establish some basic
properties
of families of $\P^1$. We prove that such a family is always locally trivial in
the \'etale
topology (prop. \ref{etaletriviality}) and we give several equivalent conditions
for the local
triviality in the Zariski topology (prop. \ref{Zariskitriviality}). Also a
cohomological interpretation
is provided in terms of the Brauer group of the base.
Surely these results are well known to the specialists but we include them here
for
the lack of an adequate bibliographical reference and also because they will
play a great role in what follows.

In section $3$, we first recall some classical basic facts about families of
hyperelliptic curves (proved in \cite{LK}): the existence of a global
hyperelliptic involution and of a
family of $\P^1$ for which the initial family of hyperelliptic curves is a
double cover, also we discuss some main properties of the Weierstrass divisor.
Then we treat the question of the existence of a global
$g^1_2$ (see the text for the precise definition). First, we give a criterion
for this existence in terms of Zariski local triviality of the underlying family
of $\P^1$
(prop. \ref{criterion}), and then we prove that such a global $g^1_2$ always
exists
if $g$ is even while for $g$ odd we give a procedure of constructing families
without
such global $g^1_2$ (theorem \ref{G_1^2}). These results were proved by
Mestrano-Ramanan (\cite{MR})
as an application of their results on Poincar\'e bundles for families of curves.
However, we believe that our approach is simpler and quite elementary.

Section $4$ deals with the moduli space $H_g$ as well as the open subset 
$H_g^0$. 
First we study the locus $H_g \setminus H_g^0$ of curves with 
extra-automorphisms 
determining the unique component of maximal dimension $g$ 
(proposition \ref{dimension}). Then 
we prove that $H_g^0$ is the smooth locus of $H_g$ except in the case $g=2$ 
where 
there is a unique singular point corresponding to the curve $y^2=x^6-x$
(proposition \ref{smooth}). After these preliminary results, we prove the two 
main theorems of this section: the determination of $\Cl(H_g)$ in theorem
\ref{picmoduli} (which turns out to be isomorphic to $\Pic(H_g^0)$) and the 
determination of $\Pic(H_g)$ in theorem \ref{picmoduli2}.

Section $5$ deals with the stack $\H$ of hyperelliptic curves. First we recall 
(including a sketch of their instructive proofs) the results of Arsie and 
Vistoli:
the description of $\H$ as a quotient stack (theorem \ref{stack}) and the 
computation of its Picard group (theorem \ref{picstack}). After that,
we provide a functorial description of a generator of the Picard group
(theorem \ref{gen}) and a description of other elements
that one can naturally consider (theorem \ref{elements}).   

In section 6 we discuss how far is the moduli functor $\mathcal{H}_g$ to be
finely represented by $H_g$.
We introduce the intermediate algebraic stack $\mathcal{D}_{2g+2}$: we describe 
it
as a quotient (theorem \ref{D-stack}), compute its Picard group and compare 
it with the Picard group of $\H$ (theorem \ref{D-pic}).
Next we prove that $\mathcal{D}_{2g+2}^0$ is indeed finely represented by
$H_g^0$ (theorem \ref{universaldivisor}) and, using this, we study how many 
families of hyperelliptic curves there can be with the same modular map 
(theorem \ref{doublecovers}). Finally we treat the existence
of a tautological family of hyperelliptic curves over an open subset 
of $H_g$ (theorem \ref{tautologicalfamily}).

The final section contains an application of the results of the preceding
section
to families of hyperelliptic curves with dominant and generically finite modular
map.
We prove that if such a family admits a global $g^1_2$ (that is the case, for
example,
if the family admits a rational section, see \cite{GV}), then the degree of the 
modular map
should be
even. This is the analog for hyperelliptic curves of a result of Caporaso
(\cite{cap03}) for families of generic smooth curves. In a forthcoming paper
(\cite{GV}), the authors
will prove an analogous result for trigonal curves and formulate a conjecture 
for $n$-gonal curves.

\emph{Acknowledgments} We are grateful to prof. A. Ragusa for organizing an
excellent summer school ``Pragmatic-2004'' held at the University of Catania,
where the two authors began their joint work on this subject.
We thank prof. L. Caporaso who suggested, during that summer school, 
an interesting research problem from which this work was originated and then
followed the progresses of this work providing useful suggestions.

\section{Families of $\P^1$}

We will call a smooth projective family of curves $p:\C\rightarrow S$ of genus
$0$ a family of
$\P^1$ (sometimes it's called a twisted $\P^1_S$ (\cite{LK}) or a conic bundle
(\cite{cil1})).
Any such family may be embedded into the projectivization of
$p_*(\omega_{\C/S}^{-1})$,
which is a vector bundle of rank $3$ on $S$.
So we obtain every family of $\P^1$ as a family of conics into a Zariski locally
trivial family of $\P^2$ (see \cite[page 12-14]{cil1}).

If the base $S$ is irreducible, the pull-back $\C_{\eta}$ of this family to the
generic point
$\eta:={\rm Spec}(k(S))\hookrightarrow S$ is a form of $\P^1_{\eta}$ (i.e. a 
variety
which becomes
isomorphic to $\P^1$ over the algebraic closure $\overline{k(S)}$).  After we
take the embedding given by the anticanonical line bundle, $\C_{\eta}$ becomes
isomorphic to a conic inside $\P^2_{\eta}$.
Recall that a conic is isomorphic to $\P^1$ if and only if it has a rational
point and surely it acquires a
rational point after a separable extension of the base field of degree $2$
(consider the field extension
given by cutting the conic with a line of $\P^2$ which intersects the conic in
two distinct points).

We want to study when $p:\C\rightarrow S$ is Zariski locally-trivial.
\begin{pro}\label{Zariskitriviality}
For the family $p:\C\rightarrow S$ over an irreducible and smooth base $S$ (with
 generic point $\eta$),
the following conditions are equivalent:
\begin{itemize}
\item[(1)] $\C$ is $S$-isomorphic to $\P(V)$ for some vector bundle sheaf $V$ on
$S$ of rank $2$.
\item[(2i)] $\C\rightarrow S$ is Zariski locally trivial.
\item[(2ii)] There exists an open non-empty $U\subset S$ such that $\C_U\cong
U\times \P^1$.
\item[(2iii)] $\C_{\eta}\cong \P^1_{\eta}$.
\item[(3ii)] The exists a rational section
$\xymatrix{
\C \ar[r]_p & S. \ar@{-->}@/_/[l]_{\sigma}
}$
\item[(3iii)] $\C_{\eta}$ has a rational point.
\item[(4i)] There exists an invertible sheaf $\mathcal{L}$ on $\C$ of vertical
degree $1$
(i.e. such that ${\rm deg}(\mathcal{L}_s)=1$ for every geometric point $s\in 
S$).
\item[(4ii)] There exists a non-empty open $U\subset S$ and an invertible sheaf
$\mathcal{L}_U$
on $\C_U$ of vertical degree $1$.
\item[(4iii)] There exists an invertible sheaf $\mathcal{L}_{\eta}$ on
$\C_{\eta}$
(i.e. defined on $\eta$) of degree $1$.
\item[(5i)] There exists an invertible sheaf $\mathcal{M}$ on $\C$ of odd
vertical degree.
\item[(5ii)] There exists a non-empty open $U\subset S$ and an invertible sheaf
$\mathcal{M}_U$
on $\C_U$ of odd vertical degree.
\item[(5iii)] There exists an invertible sheaf $\mathcal{M}_{\eta}$ on
$\C_{\eta}$ of odd degree.
\end{itemize}
\end{pro}
\begin{proof}
We will prove several implications.
\begin{itemize}
\item The implications $(\ast i)\Rightarrow (\ast ii)$ for
$\ast=2,4,5$ are evident. \item The equivalences $(\ast ii)
\Leftrightarrow (\ast iii)$ for $*=2,3,4,5$ follow from the usual
property of the generic point. \item The implications
$(1)\Rightarrow (2i)$, $(2ii)\Rightarrow (3ii)$,
$(4\star)\Rightarrow (5\star)$ (for $\star=i, ii, iii$) are
evident. $(3ii)\Rightarrow (4ii)$ follows form the fact that
$Im(\sigma)$ is the support of a divisor on $\C_U$ of vertical
degree $1$. \item $(5\star)\Rightarrow (4\star)$ (for $\star=i,
ii, iii$) follows from the fact that the relative canonical
$\omega_{\C/S}$ has vertical degree $-2$ so that, taking an
appropriate linear combination of it with $\mathcal{M}$, we obtain
an invertible sheaf $\mathcal{L}$ with vertical degree $1$. \item
$(4ii)\Rightarrow (4i)$ and $(5ii)\Rightarrow (5i)$ are true
because, thank to the smoothness of $S$ and $p$ (and hence of
$\C$), we can always extend $\mathcal{L}_U$ (or $\mathcal{M}_U$)
to an invertible sheaf on all $\C$ (simply taking the closure in
$\C$ of the Weyl=Cartier divisor in $\C_U$ corresponding to it)
and the vertical degree will remain the same since it's locally
constant and the base is connected. \item $(4i)\Rightarrow (1)$
(see \cite[prop. 3.3]{LK}):  Since the fibers of $p$ are $\P^1$,
we have that $R^1p_*(\mathcal{L})=0$ and $p_*(\mathcal{L})$ is a
locally free sheaf of rank $2$.  The natural map
$p^*(p_*(\mathcal{L}))\rightarrow \mathcal{L}$ is surjective since
its restriction to every geometric fiber is surjective. Hence it
determines an $S$-map $\Phi:\C\rightarrow \P(p_*(\mathcal{L}))$
that, being an isomorphism on the fibers, is an isomorphism.
\end{itemize}
\end{proof}

However the situation is different in the \'etale topology.

\begin{pro}\label{etaletriviality}
The family $p:\C\rightarrow S$ is locally trivial in the \'etale topology.
\end{pro}
\begin{proof}
Consider the family as a family of conics inside $\P(p_*(\omega_{\C/S}^{-1}))$.
For any point $x\in S$ we may choose a Zariski neighborhood $U$
over which $\P(p_*(\omega_{\C/S}^{-1}))$ is trivial, i.e. there is an
inclusion $\C_U\subset\P^2\times U$. Choose a line $l\subset\P^2$
that intersect the conic $\C_x$ in two different points. So there is
an \'etale double cover over some smaller Zariski neighborhood
$x\in V\subset U$ corresponding to the intersection of
$\C_V\cap (l\times V)$, over which the pull back of $\C_V$ is a
Zariski locally trivial family of $\P^1$ by \ref{Zariskitriviality} since it has
a section.
\end{proof}

There is a cohomological interpretation of this geometric picture.
The family $\C\to S$ defines a class in
$H^1_{\acute e t}(S,PGL_2(\O_S))$ by proposition \ref{etaletriviality}, and the
family is Zariski locally trivial
 if and only if it comes from $H^1_{Zar}(S,PGL_2(\O_S))$. The short exact
sequence of sheaves
$$
1\to\O_S^*\to GL_2(\O_S)\to PGL_2(\O_S)\to 1
$$
gives rise to the following two exact sequences of sheaves (for the Zariski and
the \'etale topology)
$$\xymatrix{
H^1_{Zar}(S,\O_S^*)\ar@{^{(}->}[r]\ar[d] & H^1_{Zar}(S,GL_2(\O_S)) \ar[r]\ar[d]
& H^1_{Zar}(S,PGL_2(\O_S))
 \ar[r]\ar[d]& H^2_{Zar}(S,\O_S^*)\ar[d] \\
H^1_{\acute e t}(S,\O_S^*)\ar@{^{(}->}[r]& H^1_{\acute e
t}(S,GL_2(\O_S))\ar[r] & H^1_{\acute e t}(S,PGL_2(\O_S))\ar[r] &
H^2_{\acute e t}(S,\O_S^*). }$$ It's well known that
$H^1_{Zar}(S,\O_S^*)=H^1_{\acute e t}(S,\O_S^*)=\Pic(S)$ (see
\cite[III, prop. 4.9]{Mil}) and by descent theory (see \cite[III,
sect. 4]{Mil}), there is an equality $H^1_{\acute e
t}(S,GL_2(\O_S))=H^1_{Zar}(S,GL_2(\O_S))$. Moreover for the
regular scheme $S$ the sheaf $\O_S^*$ has a flasque resolution in
the Zariski topology:
$$
1\to\O_S^*\to k(S)^*\to\oplus_{Y\in X^{(1)}}(i_Y)_*(\Z)\to 0,
$$
where $X^{(1)}$ denotes all the schematic points of codimesion 1, and $i_Y$
denotes the corresponding closed embedding. So $H^2_{Zar}(S,\O_S^*)=0$. Thus we
obtain two exact sequences:
\begin{gather}
0 \to {\rm Pic}(S) \to H^1_{Zar}(S,GL_2(\O_S)) \to H^1_{Zar}(S,PGL_2(\O_S)) \to
0 \\
0\to H^1_{Zar}(S,PGL_2(\O_S))\to H^1_{\acute e t}(S,PGL_2(\O_S))\to
H^2_{\acute e t}(S,\O_S^*). \label{Brauer}
\end{gather}
The first sequence says that every Zariski locally trivial family of $\P^1$ is
the projectivization
of a rank $2$ vector bundle (proposition \ref{Zariskitriviality}) while the
second says that
a family of $\P^1$ over $S$ defines an element in the Brauer group of $S$ which
is trivial if and
only if this family is Zariski locally trivial.

The same cohomological arguments work over the generic point
${\rm Spec}(K)$, where $K=k(S)$. Hence the exact sequence (\ref{Brauer})
can be completed in the following way
$$\xymatrix{
0\ar[r] & H^1_{Zar}(S,PGL_2(\O_S))\ar[r] & H^1_{\acute e t}(S,PGL_2(\O_S))\ar[r]
\ar[d] &
 H^2_{\acute e t}(S,\O_S^*) \ar[d]\\
& 0\ar[r] & H^1(Gal(K),PGL_2(\overline{K})) \ar[r] & {\rm Br}(K).
}$$ Since the map $H^2_{\acute e t}(S,\O_S^*)\rightarrow {\rm
Br}(k(S))$ is injective (because $S$ is smooth, see \cite[III, Ex.
2.22]{Mil}, this diagram says exactly that a family of $\P^1$
which is trivial on the generic point is Zariski locally trivial
(see proposition \ref{Zariskitriviality}).

Let us conclude this section with an example of a non-Zariski locally trivial
family of $\P^1$.
\begin{exa}\label{univconic}
Consider the universal conic ${\mathcal \C}\to S$ where $S\subset
H^0(\P^2,\O(2))$ is the open set
of all smooth conics in $\P^2$. This family is canonically embedded into
$\P^2\times S$ and is a non Zariski locally trivial family of forms
of $\P^1$, i.e. it defines a non-trivial element in ${\rm Br}(k(S))$ (see
\cite[page 16]{cil1}).
\end{exa}

\section{Generalities about families of hyperelliptic curves}

In this section, we recall first some known results about families of
hyperelliptic curves $\pi:\F\to S$,
that are projective smooth morphisms whose geometric fibers are
hyperelliptic curves of genus $g$. Recall that we assume throughout
this work that $g\geq 2$ even though many things remain true for
$g=1$ if one consider $1$-pointed elliptic curves and  family of
elliptic curves endowed with a section (see \cite{Mum} for a detailed
discussion of the elliptic case). Also recall that we work over an algebraically
closed field $k$ of characteristic different from $2$ (to avoid problems
with double covers).
\begin{teo}(\cite[Theorem 5.5]{LK})\label{hyperfamilies}
For a family $\pi:\F\to S$ of hyperelliptic curves, the following conditions
hold (and characterize
the hyperelliptic families):
\begin{itemize}
\item[(i)] $\F$ admits a global hyperelliptic involution $i$,
namely an involution over $S$ which induces the hyperelliptic
involution on every geometric fiber.

\item[(ii)] There exists a well-defined finite, surjective
$S$-morphism $h:\F\to \C$ of degree $2$, for a certain family
$p:\C\to S$ of $\P^1$ that restricts on each fiber to taking
quotient w.r.t. the hyperelliptic involution. Moreover the family
$p: \C\to S$ is uniquely determined up to $S$-isomorphisms.

\item[(iii)] The morphism $h$ may be also described as the canonical morphism
$f:\F\to \P(p_*{\omega_{\F/S}})$ whose image is isomorphic to a family
$p:\C\to S$ of $\P^1$.

\item[(iv)] There exists a faithfully flat morphism $T\to S$ and a finite
faithfully flat $T$-morphism
$\F_T\to \P^1_T=\P^1\times T$ of degree $2$.

\end{itemize}
\end{teo}

Lonsted and Kleiman studied also the Weierstrass subscheme $W_{\F/S}$ of
$\F\to S$, namely the ramification divisor of the $2:1$ of the $S$-map $h:\F\to
\C$ of theorem
\ref{hyperfamilies}(ii) endowed
with the scheme structure defined by the $0$-th Fitting ideal of
$\Omega^1_{\F/\C}$.
Note that this is isomorphic to the branch divisor $D:=h(W_{\F/S})$ on $\C$ of
the map $h$.
\begin{teo}(\cite[Prop. 6.3, Prop. 6.5, Cor. 6.8, Theo.
7.3]{LK})\label{Weierstrass}
The Weierstrass subscheme $W_{\F/S}\subset \F$ of the family of hyperelliptic
curves
$\F\to S$ satisfies the following:
\begin{itemize}
\item[(i)] $W_{\F/S}$ is the subscheme associated to an effective
Cartier divisor on $\F$ relative to $S$.
 \item[(ii)] $W_{\F/S}$ is
equal to the fixed point subscheme of $\F$ with respect to the
global hyperelliptic involution $i$. \item[(iii)] $W_{\F/S}$ is
finite and \'etale (since ${\rm char}(k)\neq 2$) over $S$ of degree
$2g+2$. \item[(iv)] If $S$ is reduced, then a section $\sigma$ of
$\pi:\F\to S$ is a Weierstrass section (i.e. $\sigma(s)$ is a
Weierstrass point of $\F_s$ for every geometric point $s\in S$) if
and only if it factors through $W_{\F/S}$.
\end{itemize}
\end{teo}
\begin{rem}\label{remark}
By the preceding results, a family $\F\to S$ of hyperelliptic
curves of genus $g$ determines a family $\C\to S$ of $\P^1$
together with a branch (Cartier) divisor $D$ on $\C$ which is
finite \'etale over $S$ of degree $2g+2$. Viceversa, by the
classical theory of cyclic covers (see for example \cite{Par} or
\cite{AV}), given the family $\C\to S$ and the divisor $D$ as
above, we can construct a double $S$-cover of $\C$ ramified
exactly over $D$ (which will be automatically a family of
hyperelliptic curves of genus $g$) if and only if the Cartier
divisor $D$ is divisible by $2$ in the Picard group of $\C$.
\end{rem}

An interesting problem for families of hyperelliptic curves is the existence of
a global $g^1_2$,
namely of an invertible sheaf $G^1_2$ on $\F$ that restricts on every fiber
of $\F\to S$ to the hyperelliptic line bundle $g^1_2$. Clearly this $G^1_2$ is
well-defined only up
to tensoring with the pull-back of line bundles coming from $S$
(see \cite[Lemma 2.1]{cil1}).
Although the uniqueness of the $g^1_2$
on every fiber of $\F\to S$ could lead to think that a $G^1_2$ always exists,
this is actually not the
case! This strange phenomenon was already observed by N. Mestrano and S. Ramanan
(\cite[section 3]{MR}) as an application of their results on Poincar\'e bundles
for families of hyperelliptic curves. Here we propose a different approach
(simpler, as we believe) that is based on the following:
\begin{pro}[Criterion for the existence of a $G^1_2$]\label{Criterion}
Let  $\pi:\F\to S$ be a family of hyperelliptic curves and let $p:\C\to S$
be a family of $\P^1$ corresponding to $\F$.
Assume $S$ is smooth and irreducible with generic point $\eta={\rm Spec}(k(S))$.
Then the following conditions are equivalent:
\begin{itemize}\label{criterion}
\item[(i)] There exists a $G^1_2$ on $\F$.
\item[(ii)] There is a non-empty open subset $U\subset S$ such that the
restriction $\F_U\to U$ admits a ${G^1_2}_{|U}$.
\item[(iii)] The hyperelliptic curve $\F_{\eta}$ admits a $g^1_2$ defined over
$k(S)$.
\item[(iv)] $p:\C\to S$ is a Zariski locally trivial family of $\P^1$.
\end{itemize}
\end{pro}
\begin{proof}
We will prove the following equivalences:
\begin{itemize}
\item[] $\underline{(i)\Leftrightarrow (ii)}$: $(i)\Rightarrow (ii)$ is clear.
Let's prove the converse.
Since $S$ and $\pi$ are smooth (and hence also $\F$), we can extend the line
bundle ${G^1_2}_{|U}$
on $\F_U$ to a line bundle $G^1_2$ on $\F$ (simply take the closure of the
Cartier=Weyl
divisor associated to it) which will have vertical degree $2$ everywhere (the
vertical degree is locally
constant and $S$ is irreducible). Now, by the semicontinuity of $h^0$
(see \cite[III.12.8]{Har}),
$h^0(\F_s,{G^1_2}_{|\F_s})\geq 2$ for every geometric point of $S$.
On the other hand for any non-zero effective divisor $E$ on an algebraic curve
$C$ there is an inequality $h^0(C,\O_C(E))\le {\rm deg}(E)$ and so in our case 
the
equality holds.
But then ${G^1_2}_{|\F_s}$ is the $g^1_2$ on the
hyperelliptic curve
 $\F_s$  (for every $s$) because this is the unique linear system of degree $2$
and dimension $1$.
\item[] $\underline{(ii)\Leftrightarrow (iii)}$: Clear from the usual property
of the generic point.
\item[] $\underline{(iii)\Leftrightarrow (iv)}$:  Consider the diagram over the
generic point $k(S)$
$$\xymatrix{
\F_{\eta} \ar[dd]_{\pi} \ar[dr]^h& \\
& \C_{\eta}\ar[ld]^p \\
k(S) & }$$ In view of proposition \ref{Zariskitriviality}, we have
to prove that $\F_{\eta}$ has a $g^1_2$ defined over $k(S)$ if and
only if $\C_{\eta}\cong \P^1_{\eta}$. Now, if $\C_{\eta}\cong
\P^1_{\eta}$ then $h^*(\O_{\P^1_{\eta}}(1))$ provides the required
$g^1_2$ on $\F_{\eta}$. Conversely, if the $g^1_2$ of $\F_{\eta}$
is defined over $k(S)$ then $V:=\pi_*(g^1_2)=H^0(g^1_2)$ is a
vector space over $k(S)$ of dimension $2$ and, by construction,
$\C_{\eta}\cong \P(V)=\P^1$.
\end{itemize}
\end{proof}

Now using this criterion, we can analyze the existence of a global $g^1_2$
(which we call
$G^1_2$) for families of hyperelliptic curves.
\begin{teo}\label{G_1^2}
Let $\F\rightarrow S$ be a family of hyperelliptic curves of genus $g$\
and let $\C\rightarrow S$ be the associated family of $\mathbb{P}^1$.\
Then the following holds:
\begin{itemize}
\item[(i)] $2 G^1_2$ (namely an invertible line bundle that
restricts to twice the $g^1_2$ on every fiber of the family) is
always defined globally on $\F\rightarrow S$. 
\item[(ii)] If $g$
is even, a $G^1_2$ is defined globally, or equivalently by the
criterion \ref{Criterion} the associated family $\C\rightarrow S$
is Zariski locally trivial. 
\item[(iii)] Viceversa, for $g$ odd,
given any family of $\mathbb{P}^1$ (maybe not Zariski locally
trivial) and a divisor \'etale and finite of degree $2g+2$ over
the base, it possible, after restricting the base to an open
subset, construct above it a family of hyperelliptic curves of
genus $g$. Hence if we start with a family of $\P^1$ non Zariski
locally trivial, the resulting family of hyperelliptic curves will
not admit a $G^1_2$ by criterion \ref{Criterion}.
\end{itemize}
\end{teo}
\begin{proof}
Keep in mind the following diagram
$$\xymatrix{
\F \ar[dd]_{\pi} \ar[dr]^h& & \\
& \C \ar[ld]^p& D \ar@{_{(}->}[l]\\
S & &
}$$
\begin{itemize}
\item[(i)] Since $\omega_{\C/S}^{-1}$ restrict to $\O(2)$ on every
fiber of $p$, the pull-back $h^*(\omega_{\C/S}^{-1})$ restricts to
$2 g^1_2$ on every fiber of  $\pi$ and hence it's the desired
$2G^1_2$. \item[(ii)] By the remark \ref{remark}, the Cartier
divisor $D$ is divisible by $2$ in the Picard group of $\C$. This
means that there exists a line bundle on $\C$ of the vertical
degree $g+1$, which is odd since $g$ is even. But then by
proposition \ref{Zariskitriviality}(5i), $\C\to S$ is Zariski
locally trivial and by the criterion \ref{criterion} there exists
a $G^1_2$ on $\F\to S$. \item[(iii)] Let $p:\C\to S$ a family of
$\P^1$ and let $D$ a divisor above it that is \'etale and finite
over $S$ of degree $2g+2$. Clearly $\O_{\C}(D)$ and
$\left(\omega_{\C/S}^{-\frac{g+1}{2}}\right)^2$ have the same
degree on the fibers so that:
$$\O_{\C}(D)\otimes \left(\omega_{\C/S}^{\frac{g+1}{2}}\right)^2
=p^*(\L)$$ for some line bundle coming from the base. Taking an
open subset $U$ of $S$ such that $\L_{|U}$ is trivial, we get that
$D$ is a square in the Picard group of $\C_U$ and therefore, by
remark \ref{remark}, we can construct the required family of
hyperelliptic curves.
\end{itemize}
\end{proof}

\section{Moduli space of hyperelliptic curves and its Picard group}

Recall that the moduli scheme $H_g$ parametrising isomorphism classes of
hyperelliptic curves is an integral subscheme of $M_g$  of
dimension $2g-1$ that can be realized as
\begin{equation}\label{descI}
H_g=(\S-\Delta)/PGL_2
\end{equation}
where $\S$ is the $(2g+2)$-th symmetric product of $\P^1$,
$\Delta$ is the closed
subset where at the least two points coincide and the action of
$PGL_2$ comes from the natural action on $\P^1$.

Equivalently, since we can identify the $(2g+2)$-th symmetric product of $\P^1$
as
the projective space $\B$ of binary forms of degree $2g+2$ in two variables, we
have
the alternative description
\begin{equation}\label{descII}
H_g=\Bsm/PGL_2
\end{equation}
where $\Bsm$ denotes the open subset of smooth binary forms (i.e. with all the
roots distinct)
and the action of $PGL_2$ is defined as $[A]\cdot [f(x)]= [f(A^{-1}x)]$.

We indicate with $H_g^0$ the open subset of $H_g$ consisting of
hyperelliptic curves with no extra-automorphisms apart from the
hyperelliptic involution. Let $\Bsm^0$ denote the preimage of $H_g^0$ inside 
$\Bsm$.

Let us remark that all the points of
$\Bsm$ are stable for the action of $PGL_2$ and with finite 
stabilizers (see \cite[prop. 4.1]{GIT}), so that the quotient
$\pi:\Bsm\to \Bsm/{PGL_2}=H_g$ is a geometric
quotient. Moreover, the action is free exactly on
$\pi^{-1}(H_g^0)=\Bsm^0$, i.e. on the forms whose corresponding 
$(2g+2)$-uples of points don't have non-trivial automorphisms.

In the next proposition, we determine the dimension of the closed subset 
$H_g^{aut}=H_g\backslash H_g^0$ of hyperelliptic
curves with extra-automorphisms and study the component of maximal 
dimension.

\begin{pro}\label{dimension}
The locus $H_g^{aut}=H_g\backslash H_g^0$ has dimension $g$ (and hence
codimension $g-1$). Moreover it has a unique irreducible component
of maximal dimension made by curves which have an extra-involution
(besides the hyperelliptic one), acting on the $2g+2$ ramification
points as a product of $g+1$ commuting transpositions.
\end{pro}
\begin{proof}
The automorphism group ${\rm Aut}(C)$ of a hyperelliptic curve $C$ always
contains the hyperelliptic involution $i$ as a central element.
Consider the group $G={\rm Aut}(C)/\langle i \rangle$. There is a
canonical inclusion inside the symmetric group $G\subset S_{2g+2}$,
since every automorphism of a hyperelliptic curve must preserve the
ramification divisor. Hence the variety $H_g^{aut}$ decomposes into
the strata
$$
H_g^{aut}=\bigcup_{{\rm primes } \: p\le 2g+2}H_g^{aut,p},
$$
where $H_g^{aut,p}$ denotes the set of hyperelliptic curves such
that there exists an element of order $p$ in the corresponding group
$G$. There is a canonical finite map $H_g^{aut,p-fixed}\to
H_g^{aut,p}$, where $H_g^{aut,p-fixed}$ is the moduli space of
isomorphism classes of pairs: a curve $C$ from $H_g^{aut,p}$ and a
fixed element $\sigma$ of order $p$ in $G$.

Since $\sigma\in G$ is induced by an automorphism of $\P^1$
preserving the ramification divisor, we see that in fact
$H_g^{aut,p-fixed}$ is the moduli space of isomorphism classes of
pairs consisting of an automorphism $\tau$ of $\P^1$ of order $p$
and a reduced effective divisor $D$ of degree $2g+2$ on $\P^1$,
stable under $\tau$. Now consider the natural quotient map
$$
\pi:\P^1=\P^1_1\stackrel{p:1}\longrightarrow\P^1_2=\P^1/\langle
\tau\rangle.
$$
The fact that $p$ is prime and the Riemann--Hurwitz formula imply
that there is only one opportunity for the ramification structure of
$\pi$: a cyclic ramification of order $p$ at two points $x_1,x_2\in
\P^1_1$. Moreover, there are three opportunities for the divisor
$D\subset \P^1_1$:
\begin{itemize}
\item[0)] $D$ contains no points among $x_1$ and $x_2$,
\item[1)] $D$ contains only one point among $x_1$ and $x_2$,
\item[2)] $D$ contains both points $x_1$ and $x_2$.
\end{itemize}
Thus we get one more stratification:
$$
H_g^{aut,p-fixed}=\bigcup_{i=0,1,2}{H_g^{aut,p-fixed,i}}
$$
according to the three cases above.

It is easy to see that in fact $H_g^{aut,p-fixed,i}$ is
parametrizing isomorphism classes of pairs, consisting of two
non-intersecting reduced effective divisors of degrees $2$ and
$(2g+2-i)/p$ on the projective line $\P^1_2$ (in this case $2g+2-i$
must be divisible by $p$). Thus, since each such configuration of
points on $\P^1$ has a finite stabilizer in the automorphism group
$PGL_2$, we get the equality
$$
\dim H_g^{aut,p-fixed,i}=2+\frac{2g+2-i}{p}-3=\frac{2g+2-i}{p}-1.
$$
Now notice that the case $p=2$ and $i=1$ is impossible because of
the divisibility condition. Further, if $p\ge 3$, or $p=2$ and
$i=2$, then
$$
\frac{2g+2-i}{p}-1\le \frac{2g+2}{3}-1\le g-1,
$$
or
$$
\frac{2g+2-2}{2}-1=g-1,
$$
respectively. So we get the inequality
$$
\dim \left(\bigcup_{{\rm primes }\: p\ge 3} H_g^{aut,p-fixed}\cup
H_g^{aut,2-fixed,2}\right)\le \max_{(p,i)\ne(2,0)} \{\dim
(H_g^{aut,p-fixed,i})\}\le g-1.
$$
If $p=2$ and $i=0$, then
$$
\dim(H_g^{aut,2-fixed,0})=g.
$$
Geometrically the condition above means that the curve $C$ has an
element $\tilde{\sigma}$ in the automorphism group ${\rm Aut}(C)$ itself
(not only in $G$). Indeed, consider the composition
$$
\varphi:C\stackrel{2:1}\longrightarrow
\P^1_1\stackrel{2:1}\longrightarrow\P^1_2.
$$
This map is a Galois map of degree 4 with Galois group $H$ generated
in ${\rm Aut}(C)$ by any preimage $\tilde{\sigma}\in {\rm Aut}(C)$ of $\sigma\in
G$ and $i$. Moreover, it is easy to see that the ramification of
$\varphi$ is formed only by pairs of double points. If $H\cong
\mathbb{Z}/4\mathbb{Z}$, then the inertia group of all the
ramification points should be the same, namely $\langle i \rangle$.
This would mean that the map $\pi\colon \P^1_1=C/\langle i
\rangle\to\P^1_2$ should be unramified, that is actually not true.
Hence $H\cong \mathbb{Z}/2\mathbb{Z}\times \mathbb{Z}/2\mathbb{Z}$,
and so $\sigma\in G$ is of order two.

Viceversa, if ${\rm Aut}(C)$ has an element $\sigma\ne i$ of order $2$
then $i=0$, otherwise $\varphi$ would have a point from
$D\cap\{x_1,x_2\}$, having ramification of order 4, contradicting
with the isomorphism $H\cong\mathbb{Z}/2\mathbb{Z}\times
\mathbb{Z}/2\mathbb{Z}$.

Note that $H_g^{aut,2-fixed,0}$ is irreducible and moreover, from
the explicit geometric description of the ramification of the
covering $C\to \P^1_2$, it follows that $\sigma\in G\subset
S_{2g+2}$ must be the product of $g+1$ commuting transpositions.
Thus we get the required statement.

\end{proof}

There exists a combinatorial proof of a weaker variant of proposition 
\ref{dimension} that we will describe now. 

\begin{pro1}
The closed subset $H_g^{aut}=H_g\backslash H_g^0$ of hyperelliptic
curves with extra-automorphisms has codimension at least 2 for $g\ge 2$, and is 
of codimension 1 for $g=2$. Moreover, in latter case the divisorial component is 
formed by hyperelliptic curves of genus 2 with an extra involution, whose action 
on the six ramification points is conjugated to $(12)(34)(56)$.
\end{pro1}
\begin{proof}
The main ingredient of the proof is the following purely combinatorial lemma.

\begin{lem}\label{combin}
Let $\rho:M\to M$ be a permutation of the finite set $M$, whose cardinality is 
at least 6. Suppose that $\rho$ has at most two fixed points, and for $|M|=6$ 
the permutation $\rho$ is not conjugated to $(12)(34)(56)$. Then there are two 
4-tuples $N_1,N_2\subset M$ such that
$$
|N_1\cap \rho(N_1)|\ne |N_2\cap \rho(N_2)|,
$$
and $|N_i\cap \rho(N_i)|<4$ for $i=1,2$ (here $|\cdot|$ denotes the cardinality 
of a set). 
\end{lem}

\begin{proof}
We treat different cases according to the cycle decomposition of $\rho$. First, 
we bound the length of cycles of $\rho$, then we bound their number, and finally 
we consider few particular cases.
\begin{itemize}
\item[{\it Case 1}] 
Suppose that there exists at least one cycle of length at least 4, i.e. there 
exists $x\in M$ such that $x_1=x$, $x_2=\rho(x)$, $x_3=\rho^2(x)$ and 
$x_4=\rho^3(x)$ are all different. Take two arbitrary elements 
$y,z\in M\backslash\{x_1,x_2,x_3,x_4\}$. One can check that 
$N_1=\{x_1,x_2,y,z\}$ and $N_2=\{x_1,x_3,y,z\}$ fit both conditions of lemma 
\ref{combin}, having $k_1=|\{x_1,y,z\}\cap\{\rho(y),\rho(z)\}|$, 
$k_2=|\{x_1,y,z\}\cap\{\rho(y),\rho(z)\}|+1\le 3$.
\item[{\it Case 2}]
Assume that there are at least four cycles. Let us take elements $x_i\in M$, 
$i=1,2,3,4$ to be in different cycles, and such that $\rho(x_1)\ne x_1$, 
$\rho(x_2)\ne x_2$. Let $l\le 2$ be equal to the number of fixed points among 
$x_3$ and $x_4$. Then for $N_1=\{x_1,x_2,x_3,x_4\}$ and 
$N_2=\{x_1,\rho(x_1),x_3,x_4\}$ there are equalities $k_1=l$ and $k_2=l+1$ or 
$l+2$, if $\rho^2(x_1)\ne x_1$ or $\rho(x_1)=x_1$, respectively. For the case, 
when $\rho(x_1)=x_1$ and $l+2=4$, i.e. when both points $x_3$ and $x_4$ are 
fixed, take $N_2=\{x_1,\rho(x_1),x_2,x_3\}$ with $k_2=3$.
\item[{\it Case 3}]
Now let us suppose that there are not more than three cycles of length at most 
3. If there are two cycles $(x_1,x_2,x_3)$ and $(y_1,y_2,y_3)$ of length 3, than 
put $N_1=\{x_1,x_2,x_3,y_1\}$, $N_2=\{x_1,x_2,y_1,y_2\}$, having $k_1=3$, 
$k_2=2$. Otherwise from the conditions on $\rho$ and $M$ we conclude that the 
lengths of cycles could be equal to $\{3,2,2\}$ or $\{3,2,1\}$. In that case 
consider the first cycle $(x_1,x_2,x_3)$ and two points $y$, $z$ from another 
two cycles such that $\rho(y)\ne y$. For $N_1=\{x_1,x_2,x_3,y\}$ and 
$N_2=\{x_1,x_2,y,z\}$ we get $k_1=3$ and $k_2=1$ or $k_2=2$.
\end{itemize}
\end{proof}

Now let $\sigma$ be an auxiliary automorphism of the hyperelliptic curve $C$ 
such that the induced permutation of $2g+2$ ramification points is not 
conjugated to $(12)(34)(56)$ for $g=2$. Then by lemma \ref{combin} there are two 
4-tuples $N_1$ and $N_2$, consisting of ramification points, such that 
$$
k_1=|N_1\cap \rho(N_1)|\ne |N_2\cap \rho(N_2)|=k_2,
$$
and $k_i<4$. Thus the point $x_C$ in $\S-\Delta=\P^6-\Delta$, corresponding to 
$C$, must lie in both closed subsets $D_{k_1}$ and $D_{k_2}$, defined in the 
following way: $D_k$ consists of points $x\in \S-\Delta$ such that in the 
corresponding $(2g+2)$-tuple of points on $\P^1$ there are two 4-tuples with the 
same double ratio and intersecting by $k$ points.

\begin{lem}\label{transvers} 
The closed subset $D_k$ is a divisor, for $k\ne 4$. Moreover, if $k_1\ne k_2$, 
and $k_i<4$ for $i=1,2$, then $D_{k_1}$ and $D_{k_2}$ intersect transversely. 
\end{lem}
\begin{proof}
Suppose that $k_1>k_2$, and $x\in D_{k_1}\cap D_{k_2}$. For a $(2g+2)$-tuple 
$M_x$ of points on $\P^1$, corresponding to $x$, there exists at least one pair 
of 4-tuples $(N_1,N_2)$ with the same double ratio, such that $|N_1\cap 
N_2|=k_1$. Consider also all pairs $(L^i_1,L^i_2)$ of 4-tuples in $M_x$ with the 
same double ratio, such that $|L^i_1\cap L^i_2|=k_2$. For any $i$ there exists a 
point $z_i\in M_x$ such that $z_i\in L^i_1\cup L^i_2$, but $z_i\notin N_1\cup 
N_2$, since $|L^i_1\cup L^i_2|=8-k_2>8-k_1=|N_1\cup N_2|$. So for each point 
$x\in D_{k_1}\cap D_{k_2}$ and for each Zariski neighborhood $x\in U\subset 
D_{k_1}$ there exists a point $y\in U$ not belonging to $D_{k_2}$: just slightly 
move all the points $z_i$ in an independent way. This provides the 
transversality of the intersection $D_{k_1}\cap D_{k_2}$.
\end{proof}

From lemma \ref{transvers} we get that the set of hyperelliptic curves with 
auxiliary automorphisms is contained inside the closed subset
$$
\bigcup_{0\le k_1<k_2< 4}(D_{k_1}\cap D_{k_2}),
$$
which is of codimension 2, if $g\ge 3$, or $\rho$ is not conjugated to 
$(12)(34)(56)$ for $g=2$. 

Now consider the case, when $g=2$ and $\rho$ is conjugated to $(12)(34)(56)$. 
The dimension of the moduli space of such hyperelliptic curves (i.e. 
hyperelliptic curves of genus 2 having an auxiliary automorphism with the action 
of type $(12)(34)(56)$ on the ramification points) is equal to 5=3+2: six 
ramification points are uniquely defined by three of them and also the 
involution in $PGL_2$ that is a two-dimensional space, since the involutions are 
parameterized by a couple of their fixed points. Thus we get the desired 
statement.
\end{proof}

\begin{rem}
For the moduli space $M_g$ of curves of genus $g$ (with $g\geq 3$),  the 
locus $M_g^{aut}$ of curves with non-trivial automorphisms is a closed subset 
of dimension $2g-1$ and it has a unique irreducible component of maximal 
dimension made by hyperelliptic curves.
\end{rem}

Note that $H_g$ is a normal variety since it is the quotient of a normal variety 
by the action of a group. We determine its smooth locus. 

\begin{pro}\label{smooth}
If $g\geq 3$, then the smooth locus of $H_g$ is $H_g^0$. On the other hand, 
assuming ${\rm char}(k)\neq 5$, $H_2$ has a unique singular point corresponding 
to the hyperelliptic curve  $y^2=x^6-x$. 
\end{pro}
\begin{proof}
We will use a smoothness criterion of K. Lonsted (\cite{Lon1})
for quotients of smooth varieties by finite groups, that we now briefly recall.

Let $X$ be a smooth $k$-variety,
$\Gamma\subset {\rm Aut}_k(X)$ a finite subgroup, and set $Y=X/\Gamma$. Let $P$
be a point in $X$ with image $Q$ in $Y$. Let's denote with $\Gamma(P)$ the 
inertia 
group at $P$ and let $\Gamma'(P)$ be the subgroup of $\Gamma(P)$ generated 
by all the pseudoreflections, that is, the elements of $\Gamma(P)$ that leave a 
hypersurface 
through $P$ pointwise fixed. Then one has the following: 
\begin{itemize}
\item[(i)] If $Q$ is a smooth point, then $\Gamma(P)=\Gamma'(P)$.
\item[(ii)] Conversely, if $\Gamma(P)=\Gamma'(P)$ and the order of $\Gamma(P)$ 
is
prime with ${\rm char}(k)$, then $Q$ is a smooth.  
\end{itemize}

Note that $H_g$ can be realized as a quotient of a smooth variety by a finite 
group in the 
following way. Given a $(2g+2)$-uple of ordered distinct points of $\P^1$, 
acting 
with an element of $PGL_2$ we can assume that the first three points of it 
are $0, \infty, 1$. Hence
\begin{equation}\label{descIII}
H_g=((\P^1-\{0, \infty, 1\})^{2g-1}-\Delta)/S_{2g+2}
\end{equation}
where $\Delta$ is the locus where at least two points coincide and the action of 
an element $\sigma\in S_{2g+2}$ on an ordered $(2g-1)$-tuple $\{x_1,\cdots, 
x_{2g-1}\}$ 
is obtained first letting $\sigma$ act in the natural way on the $(2g+2)$-tuple 
$\{0, \infty, 1, x_1, \cdots, x_{2g-1}\}$
and then applying the element of $PGL_2$ that sends the first three elements 
into 
$\{0, \infty, 1\}$ and taking the remaining $(2g-1)$ points.

Apply the preceding smoothness criterion with $X= (\P^1-\{0, \infty, 1\})^{2g-
1}-\Delta$
and $\Gamma=S_{2g+2}$. 
If $g\geq 3$, proposition \ref{dimension} implies that there aren't
non-trivial pseudoreflections. In fact such a non-trivial pseudoreflection would 
imply
the existence of a hypersurface on $X$ made by points having non-trivial 
stabilizer and, 
passing to the quotient, this would give a codimension $1$ locus of 
hyperelliptic curves 
with extra-automorphisms contradicting proposition \ref{dimension}. Hence, by 
the criterion,
a point on the quotient $H_g$ is non-singular if and only if it comes from a 
point above 
with trivial stabilizer, hence if and only if it belongs to $H_g^0$. 

If $g=2$, this argument fails because in that case the elements of $S_6$ 
conjugated to 
$(12)(34)(56)$ are pseudoreflections (and by proposition \ref{dimension}, these 
are the 
only ones). In this case we can use, instead, an explicit description
of Igusa (see \cite{Igu}) who showed that (under the hypothesis 
${\rm char}(k)\neq 5$):
\begin{equation}\label{Igusa}
M_2=H_2=
{\rm Spec}(k[z_1,z_2,z_3])/{\langle\zeta_5\rangle}=\mathbb{A}_k^3/(\Z/5\Z)
\end{equation}
where the action of the $5$-th root of unity $\zeta_5$ is given by
$z_i\mapsto \zeta_5^i z_i$, $i=1, 2, 3$ and the origin corresponds to the 
hyperelliptic curve defined by the equation $y^2=x^6-x$. It is well-known that 
the origin in $\mathbb{A}_k^3$ is mapped to the singular point on $H_2$. Also we 
could get it applying the smoothness criterion to this quotient, since in this 
case there aren't pseudoreflections, it follows that the only singularity of the 
quotient is the point corresponding to the curve $y^2=x^6-x$.  
\end{proof}

\begin{rem}
Compare this result with the determination of the smooth locus of $M_g$ for 
$g\geq 3$.
In this case it holds that $M_g^{smooth}=M_g^0$ for $g\geq 4$ (i.e. exactly 
when the locus of curves with automorphisms has codimension greater than $1$)
while $M_3^{smooth}=M_3^0\cup H_3^0$ (see \cite{Rau} for an analytic proof over 
the 
complex numbers, \cite{Pop} for an algebraic proof in the case $g\geq 4$, 
\cite{Oort}
for an algebraic proof in the case $g=3$ and finally \cite{Lon2} for an 
algebraic 
unified treatment of the cases $g=3$ and $g\geq 4$ based on its smoothness 
criterion 
\cite{Lon1}). 
\end{rem}

We want now to compute the Picard groups (i.e. the group of Cartier divisors 
modulo 
linear equivalence) and the divisor class groups (i.e. the group of Weyl 
divisors 
modulo linear equivalence) of $H_g$ and of $H_g^0$, away from 
some bad characteristic of the base field $k$. 
Note that since $H_g$ is a normal variety we have an inclusion 
$\Pic(H_g)\hookrightarrow \Cl(H_g)$; 
on the other hand,
$H_g^0$ is smooth (see proposition \ref{smooth}) and hence 
$\Pic(H_g^0)=\Cl(H_g^0)$. 

\begin{teo}\label{picmoduli}
Suppose that ${\rm char}(k)$ doesn't divide $2g+2$.
The Picard group of $H_g^0$ is equal to
$$\Pic(H_g^0)=\left\{
\begin{aligned}
&\mathbb{Z}/(4g+2)\mathbb{Z} \text{ if } g \geq 3 \\
&\mathbb{Z}/5\mathbb{Z} \text{ if } g=2.\\
\end{aligned}
\right.$$
Moreover, under the additional hypothesis that ${\rm char}(k)\neq 5$ if $g=2$,
the natural restriction map $\Cl(H_g)\to \Cl(H_g^0)\cong \Pic(H_g^0)$ is an
isomorphism.
\end{teo}
\begin{proof}
We will use the theory of \emph{equivariant Picard group} (see
\cite[I,3]{GIT} and also \cite{EG}) whose definition we now briefly recall.
Given an action of an algebraic group $G$ on a algebraic variety
$X$, $\sigma:G \times X\to X$, the equivariant Picard group
$\Pic^G(X)$ is defined as
$$\Pic^G(X)=\{(\L,\phi): \L\in \Pic(X), \phi \mbox{ is a }G-
\mbox{linearization}\}_{/\cong}$$
where a $G$-linearization $\phi$ of a line bundle $\L$ is an
isomorphism $\phi:\sigma^*(\L)\stackrel{\cong}{\rightarrow}
p_2^*(\L)$ ($p_2$ is the projection $G\times X\to X$) satisfying
the obvious cocycle condition (see \cite[pag. 30]{GIT}).
We will apply this
in our case with $G=PGL_2$ and $X=\Bsm$ (see \ref{descII}).

In this case, since there aren't non-trivial homomorphisms
$PGL_2\to \mathbb{G}_m$, we have an injection
$\Pic^{PGL_2}(\Bsm) \hookrightarrow \Pic(\Bsm)$
(see \cite[prop. 1.4]{GIT}).

Moreover, since $\Delta=\B-\Bsm$ is an irreducible hypersurface $\Delta$ of
degree $4g+2$ (lemma \ref{Delta}), from the
exact sequence (see \cite[II, 6.5]{Har})
$$\Z\cdot [\Delta]\to \Pic(\B)\to \Pic(\Bsm)\to 0$$
we get that $\Pic(\Bsm)=\Z/(4g+2)\Z$ generated by the
hyperplane section $\O(1):=\O_{\B}(1)_{|\Bsm}$.

$\underline{CLAIM}:$ $\O(1)$ admits a $PGL_2$-linearization.\\
In fact since the action of $\sigma:PGL_2\times \B\to \B$ is
linear in $\B$ and of degree $2g+2$ in $PGL_2$, we have that
$$\sigma^*(\O_{\B}(1))=p_1^*(\O_{PGL_2}(2g+2))\otimes p_2^*(\O_{\B}(1)).$$
Moreover since $PGL_2=\P^4-\{{\rm det}=0\}$ and ${\rm det}$ is of
degree $2$, $\Pic(PGL_2)=\Z/2\Z$ and hence
$\O_{PGL_2}(2g+2)\cong \O_{PGL_2}$. From this, it follows that
$\sigma^*(\O_{\B}(1))\stackrel{\cong}{\rightarrow}
p_2^*(\O_{\B}(1))$ and hence the claim. So we reached the
conclusion that
\begin{equation}\label{4.2}
\Pic^{PGL_2}(\Bsm)=\Z/(4g+2)\Z
\end{equation}
generated by $\O_{\B}(1)$ (for every $g\geq 2$).

The last statement has another explanation: if an algebraic group
$G$ acts on the projective space $\P^n=\P(V)$, then the sheaf
$\O_{\P^n}(1)$ admits a $G$-linearization if and only if the initial
action is induced from a representation of $G$ in the vector space
$V$. It follows from the inclusion of the tautological bundle
$\O_{\P^n}(-1)$ into the product $\P^n\times V$ and the diagonal
action of $G$ on $\P^n\times V$. In our case $SL_2$ does act on the
vector space of binary forms of degree $2g+2$ in two variables by
the same formula as $PGL_2$ on the projective space. Moreover,
$\{\pm 1\}={\rm Ker}(SL_2\to PGL_2)$ acts trivially on the binary forms of
even degree, so $PGL_2$ also acts on this vector space, and hence on
$\O_{\B}(1)$. Explicitly the action of a class $[A]$ of $PGL_2$ on a binary 
form $f(x)$ is given by: $[A]\cdot f(x)={\rm det}(A)^{g+1}f(A^{-1}\cdot x)$.   

Now we are going to relate this equivariant Picard group with the
divisor class group of the quotient variety $H_g=\Bsm/{PGL_2}$
(note that a priori the equivariant Picard group is the Picard
group of the quotient stack $[\Bsm/{PGL_2}]$ (see
\cite[prop. 18]{EG})). 

Using the theory of descent, one can show that if the action is
free and the quotient is a geometric quotient then the equivariant
Picard group is the Picard group of the quotient variety (see
\cite[pag. 32]{GIT}), so that in our case:
$$
\Pic(H_g^0)=\Pic^{PGL_2}(\Bsm^0).
$$

Now, if $g\geq 3$, the proposition \ref{dimension} says that
$H_g-H_g^0$ has codimension greater or equal to $2$. From this, it
follows that $\Cl(H_g)\stackrel{\cong}{\to}\Cl(H_g^0)$ (see \cite[II.6.5]{Har}), 
and
$\Pic^{PGL_2}(\Bsm)=\Pic^{PGL_2}(\Bsm^0)$ (see
\cite[sect. 2.4, lem. 2]{EG}) so that:
\begin{equation}\label{4.4}
\Cl(H_g)\stackrel{\cong}{\to}\Pic(H_g^0)=\Pic^{PGL_2}(\Bsm)
\end{equation}
which together with \ref{4.2} gives the conclusion.

If $g=2$, this argument fails because in this case $H_2-H_2^0$ contains
a divisor and removing it affects the divisor class group.
We will compute $\Cl(H_2)$ and $\Pic(H_2^0)$ in two different ways
obtaining that they are both isomorphic to $\Z/5\Z$.

First of all, to compute $\Cl(H_2)$ we use the explicit description
of Igusa under the hypothesis ${\rm char}(k)\neq 5$ (see formula \ref{Igusa}). 
Since the action of $<\zeta_5>$ is free outside the point $C_0:=\{y^2=x^6-x\}$
(that has codimension $3$), the same reasoning as before gives
$$\Cl(H_2)=\Cl(H_2-[C_0])=\Pic(H_2-[C_0])=\Pic^{\Z/5\Z}(\mathbb{A}_k^3-
0)=\Pic^{\Z/5\Z}(\mathbb{A}_k^3)\cong \Z/5\Z.$$

Next, let $\widetilde{D}$ be the unique irreducible component of codimension
$1$ of $H_2-H_2^0$ (see proposition \ref{dimension}) and let $D$ its inverse 
image
in $\mathbb{B}_{sm}(2,6)=\mathbb{B}(2,6)-\Delta$ and $\overline{D}$ its closure
in
$\mathbb{B}(2,6)$.
The same reasoning as before shows that
$$\Pic(H_2^0)=\Pic^{PGL_2}(\mathbb{B}(2,6)-(\Delta\cup \overline{D})).$$
By lemma \ref{Delta}, $\Delta$ is an irreducible hypersurface of degree $10$ in
$\mathbb{B}(2,6)\cong\P^6$ and, by lemma \ref{g=2}, $\overline{D}$ is
an irreducible hypersurface of degree $15$, so that
$\Pic(\mathbb{B}(2,6)-(\Delta\cup \overline{D}))\cong \Z/5\Z$ (use the usual
exact sequence of \cite[II.6.5]{Har}).
Moreover from the claim above it follows that
$\Pic^{PGL_2}(\mathbb{B}(2,6)-(\Delta\cup \overline{D}))=
\Pic(\mathbb{B}(2,6)-(\Delta\cup \overline{D}))$ and hence
the desired conclusion.
\end{proof}

\begin{lem}\label{Delta}
Suppose that ${\rm char}(k)$ doesn't divide $2g+2$. The closed subset
$\Delta=\B-\Bsm$ (given by the vanishing of the discriminant) is an
irreducible hypersurface of degree $4g+2$ in $\B=\P^{2g+2}$.
\end{lem}
\begin{proof}
Let's consider the polynomial $f$ of degree $n:=2g+2$ associated to
a binary form. Recall that the discriminant $\Delta(f)$ is the
resultant $R(f,f')$ of the polynomial with its derivative divided by
the leading coefficient (see \cite[pag. 104]{GKZ}). The resultant
$R(f,f')$ is the determinant of a square matrix of size $n+n-1=2n-1$
whose entries are the coefficients of our polynomial and hence it
will be a homogeneous polynomial in these coefficients of degree
$2n-1$. It follows that the discriminant will be homogeneous of
degree $2n-2$ which in our case gives $4g+2$ (for another proof see
\cite{Ran}).

The irreducibility of the discriminant polynomial (under the
hypothesis that ${\rm char}(k)$ doesn't divide $2g+2$) is proved in
\cite[pag. 658-659]{AV}.
\end{proof}

\begin{lem}\label{g=2}
Let $D$ be the unique irreducible component of codimension $1$ of
$\BSsm-\BSsm^0$ (see proposition \ref{dimension}) and let $\overline{D}$
be its closure in $\BS$. Then $\overline{D}$ is an irreducible
hypersurface in $\BS=\P^6$ of degree $15$.
\end{lem}
\begin{proof}
Let us consider the map
$\pi:(\P^1)^6-\Delta\stackrel{S_6}{\longrightarrow}
{\rm Sym}^6(\P^1)-\Delta$, where $\Delta$ indicates in both spaces
(with an abuse of notation) the locus of $6$-tuples of points with
at least $2$ coincident points.

We want to decompose the divisor $\pi^{-1}(D)$ in $(\P^1)^6-\Delta$
or, more precisely, its closure $\overline{\pi^{-1}(D)}$ in
$(\P^1)^6$. By proposition \ref{dimension}, an element of $\pi^{-1}(D)$ is
a $6$-tuple of distinct ordered points of $\P^1$ that has an
automorphism of order $2$, whose action on these six points is
conjugated to $(12)(34)(56)$, or in other words such that there
exists an element $A\in PGL_2$, inducing such permutation $\sigma$
of the 6-tuple. So we obtain a decomposition
\begin{equation}\label{decomposition}
\pi^{-1}(\overline{D})=\bigcup_{\sigma\sim (12)(34)(56)}
\overline{D}_{\sigma}
\end{equation}
where the union is taken over the $15$ elements of $S_6$ conjugated
to $(12)(34)(56)$, and for each of them $\overline{D}_{\sigma}$ is
an hypersurface.

Now we will compute the class of $\overline{D}_{\sigma}$ in the
Picard group $\Pic((\P^1)^6)\cong(\mathbb{Z})^6$ (without loss of
generality we can consider $D_{(12)(34)(56)}$). Take a line
$l=\{P_1\}\times\ldots\times\{P_5\}\times\P^1$ in $(\P^1)^6$ for
general points $P_i\in\P^1$. Let $P=(P_1,\ldots,P_5,P_6)\in l\cap
D_{(12)(34)(56)}$, and let $A\in PGL_2$ be an automorphism, inducing
the corresponding permutation of $P_i$. We have the following
conditions on $A$:
$$\begin{sis}
A(P_1)=P_2,\\
A(P_3)=P_4,\\
A^2=1.\\
\end{sis}
$$
The point $P_6=A(P_5)$ is uniquely determined by $A$, so we want to
understand how many $A$ are satisfying the conditions above.

Choose the homogenous coordinates of $P_1$, $P_2$, $P_3$ and $P_4$
to be equal to $[(1:0)]$, $[(1:1)]$, $[(0:1)]$ and $[(c:d)]$
respectively (with $c\ne 0$, $d\ne 0$, $c\ne d$). Then, due to the
first two conditions, the matrix $A$ should be equal in this basis
to
$$
A=\left[\begin{pmatrix}
1&\lambda c\\
1&\lambda d\\
\end{pmatrix}\right]
$$
for some nonzero $\lambda$. The last condition $A^2=1$ gives
$\lambda=-1/d$. Besides, since the $P_i$ are general, $A(P_5)\ne
P_1, P_2, P_3, P_4, P_5$, so $l\cap D_{(12)(34)(56)}$ consists of
one point, that is a transversal intersection.

Moreover, the five points from the intersection $l\cap \Delta$
cannot lie on $\overline{D}_{(12)(34)(56)}$: if a point
$Q=(P_1,\ldots,P_5,Q_6)\in l\cap \Delta $ is a limit of points
$Q^t\in \overline{D}_{(12)(34)(56)}$ then at each moment $t$ the
point $Q^t_6\in \P^1$ is uniquely algebraically determined by
$Q^t_5$ and $A^t\in PGL_2$, that is uniquely algebraically
determined by $(Q_1^t,Q_2^t,Q_3^t,Q_4^t)$. Hence $Q_6$ must be equal
to $P_6$, so $l\cap (\overline{D}_{(12)(34)(56)}- D_{(12)(34)(56)})$
is empty.

Now due to the symmetry of $\overline{D}_{(12)(34)(56)}$ the same is
true for all other ``coordinate'' lines in $(\P^1)^6$, and so the
class of $\overline{D}_{(12)(34)(56)}$ in $\Pic((\P^1)^6)$ is equal
to $(1,1,1,1,1,1)$. Thus, combining this result with the
decomposition (\ref{decomposition}) and comparing it with the fact
that $\pi^{-1}(\O_{\P^6}(1))$ is also of type $(1,1,1,1,1,1)$, we
obtain that the degree of $\overline{D}$ is equal to 15.

\end{proof}

\begin{teo}\label{picmoduli2}
Suppose that ${\rm char}(k)$ doesn't divide $2g+2$ or $2g+1$.
Then 
$$\Pic(H_g)=0.$$
\end{teo}
\begin{proof}
Consider the following natural maps (see theorem \ref{picmoduli}):
$$\Pic(H_g)\to\Pic^{PGL_2}(\Bsm)\twoheadrightarrow\Pic^{PGL_2}(\Bsm^0)
\cong\Cl(H_g).$$
Since $H_g$ is normal, the composition of the two maps is an injection, hence
also the first map is an injection. 

Recall (see formula \ref{4.2}) that $\Pic^{PGL_2}(\Bsm)$ is a cyclic group of 
order
$4g+2$ generated by the tautological line bundle $\O_{\Bsm}(-1)$ with its 
natural $PGL_2$-linearization that comes from its embedding inside 
$\Bsm\times \Asm$, where $PGL_2$ acts diagonally.

We want to see which $PGL_2$-linearized line bundles $L$ on  $\Bsm$ come from 
line bundles 
on $H_g$.  Clearly a necessary condition is that for each point 
$x\in\Bsm$ its stabilizer ${\rm Stab}_x\subset PGL_2$ is acting trivially on the 
fiber $L|_x$. 

Consider first the binary form $f_1:=X^{2g+1}Y-Y^{2g+2}$ (which is in $\Bsm$ 
since 
${\rm char}(k)$ doesn't divide $2g+1$). Its stabilizer is the cyclic group of 
order 
$2g+1$:
$${\rm Stab}_{f_1}=C_{2g+1}=\left\langle \left[
\begin{array}{cc}
\zeta_{2g+1}& 0\\
0& 1\\
\end{array}\right] \right\rangle$$
where $\zeta_{2g+1}$ is a primitive $(2g+1)$-root of unity.

The fiber of the line bundle $\O_{\Bsm}(-1)$ above $f_1$ is the $1$-dimensional 
vector 
space of all scalar multiples of $f_1$ inside $\Asm$:
$$\O_{\Bsm}(-1)_{f_1}=\{\lambda\cdot (X^{2g+1}Y-Y^{2g+2}) : \lambda\in k\}.$$
Recall (from the proof of theorem \ref{picmoduli}) that $PGL_2$ acts on $\Asm$ 
by 
the formula: $[A]\cdot f(x)={\rm det}(A)^{g+1}f(A^{-1}\cdot x)$. So the 
generator 
of the stabilizer group acts on the fiber as multiplication by 
$\zeta_{2g+1}^{g+1}=
\zeta_{2g+1}^{-g}$. Hence only the multiples of $\O_{\Bsm}(2g+1)$ can come from 
line 
bundles on $H_g$.

Next consider the binary form $f_2:=X^{2g+2}-Y^{2g+2}$ (which is in $\Bsm$ since 
${\rm char}(k)$ doesn't divide $2g+2$). Its stabilizer is the diedral group of 
order 
$4g+4$:
$${\rm Stab}_{f_2}=D_{2g+2}=\left\langle 
\left[
\begin{array}{cc}
\zeta_{2g+2}& 0\\
0& 1\\
\end{array}\right], 
\left[
\begin{array}{cc}
0& 1\\
1& 0\\
\end{array}\right]
 \right\rangle$$
where $\zeta_{2g+2}$ is a primitive $(2g+2)$-root of unity.

The fiber of the line bundle $\O_{\Bsm}(-1)$ above $f_2$ is:
$$\O_{\Bsm}(-1)_{f_2}=\{\lambda\cdot (X^{2g+2}-Y^{2g+2}) : \lambda\in k\}.$$
The two generators of the stabilizer group act respectively as multiplication
by $-1$ and $(-1)^g$. Hence only the multiples of $\O_{\Bsm}(2)$ can come 
from line bundles on $H_g$.

Putting together these two conditions plus the fact that $\O_{\Bsm}(4g+2)=0$ in 
$\Pic^{PGL_2}(\Bsm)$,  one concludes that $\Pic(H_g)=0$. 
\end{proof}

\section{Stack of hyperelliptic curves and its Picard group}

Recall that the moduli functor $\mathcal{H}_g$ of hyperelliptic curves 
is the contravariant functor
$$\mathcal{H}_g : Sch_{/k}\rightarrow Set$$
which associates to every $k$-scheme $S$ the set
$$\mathcal{H}_g(S)=\left\{
\F\rightarrow S \text{ family of hyperelliptic smooth curves of genus } g
\right\}_{/\cong}.$$

By the results of Lonsted-Kleiman (see theorem \ref{hyperfamilies}
and \ref{Weierstrass}), a family $\pi:\F\to S$ of hyperelliptic curves is 
a double cover of a family $p:\C\to S$ of $\P^1$, namely we have the following 
situation 
$$\xymatrix{
W \subset \mathcal{F} \ar[dr]^{f}_{2:1} \ar@<1ex>[dd]_{\pi} & \\
& \C\supset D \ar[dl]_{\P^1}^p \\
S &
}$$
where the branch divisor $D$ and the ramification divisor $W$ (the Weierstrass 
subscheme) are relative Cartier divisor finite and \'etale of degree $2g+2$ 
over the $S$. By the classical 
theory of double covers, the divisor $D$ is divisible by $2$ in the Picard group 
of 
$\C$, namely there exists an invertible sheaf $\L$  in ${\rm Pic}(\C)$ such that 
\begin{equation}\label{f1}
(\L^{-1})^{\otimes 2}=\O_{\C}(D).
\end{equation}
This invertible sheaf satisfies the following two relations
\begin{equation}\label{f2}
f^*(\L^{-1})=\O_{\F}(W)
\end{equation}
\begin{equation}\label{f3}
f_*(\O_{\F})=\O_{\C}\oplus \L.
\end{equation}
Moreover the Hurwitz formula gives
\begin{equation}\label{f4}
\omega_{\F/S}=f^*(\omega_{\C/S})\otimes \O_{\F}(W).
\end{equation} 

In view of these results, one can prove (see \cite[section 3]{AV}) that the 
functor
$\H$ is isomorphic to the functor $\H'$  which 
associates to a $k$-scheme $S$ the set
$$
\H'(S)=\{\C\stackrel{p}{\to} S, \L, \L^{\otimes 2}\stackrel{i}{\hookrightarrow} 
\O_{\C}\}
$$
where $p:\C\to S$ is a family of $\P^1$, $\L$ is an invertible sheaf on $\C$ 
that restricts to an invertible sheaf of degree $-g-1$ on any geometric fiber
and $i:\L^{\otimes 2}\hookrightarrow \O_{\C}$ is an injective map of line 
bundles
that remains injective on any geometric fiber and such that the image of $i$ is
the sheaf of ideals of a relative Cartier divisor finite and \'etale over $S$ 
(it's the branch divisor $D$ in the description above).

In \cite{AV}, Arsie and Vistoli proved that $\H'\cong \H$ is a Deligne-Mumford 
algebraic stack and describe it as a quotient stack (more generally, they 
consider
stacks of cyclic covers of projective spaces). This explicit description 
allows them to compute the Picard group of $\mathcal{H}_g$ 
(in the sense of Mumford \cite{Mum}). We are going to recall their results here.

Consider, inside the affine space $\A$ of linear forms in two variables of
degree $2g+2$, the open subset $\Asm$ of smooth linear forms 
(i.e. forms having distinct roots)
and an action of $GL_2$ by: $A\cdot f(x)= f(A^{-1}\cdot x)$.  Let us remark that
the projective space $\B$ is just the projectivization of $\A$ (the same is true
for $\Bsm$ and $\Asm$). Clearly the subgroup
$\mu_{g+1}$, embedded diagonally in $GL_2$, acts trivially on $\Asm$.
The result is the following:
\begin{teo}(Arsie-Vistoli, \cite[theo. 4.1]{AV})\label{stack}
The stack of hyperelliptic curves of genus $g$ can be realized as
$$\mathcal{H}_g=[\Asm/(GL_2/\mu_{g+1})]$$
with action given by $[A]\cdot f(x)=f(A^{-1}\cdot x)$.
\end{teo}
\begin{proof}(Sketch)
Consider the auxiliary functor $\widetilde{\H}$ which associates to every $k$-
schemes $S$
the set
$$
\widetilde{\H}(S)=\{\C\stackrel{p}{\to} S, \L, \L^{\otimes 
2}\stackrel{i}{\hookrightarrow} \O_{\C}, \phi:(\C,\L)\cong (\P^1_S,\O_{\P^1_S}
(-g-1))\}
$$
where $p:\C\to S$, $\L$ and $i:\L\hookrightarrow \O_{\C}$ are as before and 
the isomorphism $\phi$ consists of an isomorphisms of $S$-schemes 
$\phi_0:\C\cong \P^1_S$ plus an isomorphism of invertible sheaves
$\phi_1:\L\cong \phi_0^*\O_{\P^1_S}(-g-1)$. Forgetting the isomorphism $\phi$, 
one 
gets a natural transformation of functors $\widetilde{\H}\to \H'\cong \H$. 

This rigidified functor $\widetilde{\H}$ is isomorphic to $\Asm$ (thought as the 
functor
${\rm Hom}(-,\Asm)$). In fact for any object in $\widetilde{\H}(S)$ the 
isomorphism $\phi_1$ (precisely, its ``tensor square'') and the inclusion 
$i:\L^{\otimes 2}\hookrightarrow \O_{\C}$ provide a canonical inclusion 
$\O_{\P^1_S}(-2g-2)\hookrightarrow \O_{\P^1_S}$. This morphism of sheaves 
corresponds to a section of $\O_{\P^1_S}(2g+2)$, that is smooth on any geometric 
fiber 
of $\P^1_S\to S$, and so defines an element of $\Asm(S)$. The inverse functor is 
obtained 
by sending an element of $\Asm(S)$, thought as a homomorphism 
$f:\O_{\P^1_S}(-2g-2)\to \O_{\P^1_S}$, into the object
$\{\P^1_S\to S, \O_{\P^1_S}(-g-1), 
f:\O_{\P^1_S}(-g-1)^{\otimes 2}\to\O_{\P^1_S},
{\rm id}:(\P^1_S, \O_{\P^1_S}(-g-1))\to
(\P^1_S,\O_{\P^1_S}(-g-1))\}$
of $\widetilde{\H}(S)$.  

Next consider the group ${\rm Aut}(\P^1,\O(-g-1))$ consisting of automorphisms 
of $\P^1$ with a linearization of the sheaf $\O(-g-1)$. This group is 
canonically isomorphic to $GL_2/\mu_{g+1}$, and the corresponding group sheaf 
$\underline{{\rm Aut}}(\P^1,\O(-g-1))$ 
acts naturally on $\widetilde{\H}$ by composition with the isomorphism $\phi$. 
One can check that the corresponding action of $GL_2/\mu_{g+1}$ on $\Asm$ is the 
one given in the statement.

Finally, descent theory implies that the forgetful morphism $\widetilde{\H}\to 
\H$ 
makes $\widetilde{\H}$ into a principal bundle over the stack $\H$ respect to
the group sheaf  $\underline{{\rm Aut}}(\P^1,\O(-g-1))$. 
From this, one gets the representation of
$\H\cong \H'$ as the quotient stack $[\Asm/(GL_2/\mu_{g+1})]$. 
\end{proof}

Note that also the bigger group $\mu_{2g+2}$ acts trivially on $\Asm$ 
the stack $[\Asm/(GL_2/{\mu_{2g+2}})]$ is not isomorphic to
$[\Asm/(GL_2/{\mu_{g+1}})]$ (see remark \ref{D-stack2}).

One can give a more explicit description of the quotient group 
appearing in the preceding theorem as in the following: 
\begin{lem}\label{quotientgroup}
For the group $GL_2/ \mu_{g+1}$ it holds:
\begin{itemize}
\item[(i)] If $g$ is even then the homomorphism of algebraic groups
$$GL_2/\mu_{g+1}\to GL_2$$
given by
$$[A]\mapsto {\rm det}(A)^{\frac{g}{2}}A$$
is an isomorphism. The group of characters of $GL_2/\mu_{g+1}$ is 
isomorphic to $\Z$ and is generated by ${\rm det}^{g+1}$.
\item[(ii)] If $g$ is odd then the homomorphism of algebraic groups
$$GL_2/\mu_{g+1}\to \mathbb{G}_m\times PGL_2$$
given by
$$[A]\mapsto ({\rm det}(A)^{\frac{g+1}{2}},[A])$$
is an isomorphism. The group of characters of $GL_2/\mu_{g+1}$ is 
isomorphic to $\Z$ and is generated by ${\rm det}^{\frac{g+1}{2}}$.
\end{itemize}
\end{lem}
\begin{proof}
\begin{itemize}
\item[(i)] An inverse is given by the homomorhism
$$A\mapsto [{\rm det}(A)^{\frac{-g}{2(g+1)}}A].$$
The second assertion follows from the fact that the group of characters of 
$GL_2$ is
isomorphic to $\Z$ and is generated by ${\rm det}$.
\item[(ii)] An inverse is given by the homomorphism
$$(\alpha, [A])\mapsto [\alpha^{\frac{1}{g+1}}{\rm det}(A)^{-1/ 2}A].$$ 
The second assertion follows from the fact that the group of characters of 
$\mathbb{G}_m\times PGL_2$ is isomorphic to $\Z$ and is generated by the 
projection onto the first factor.
\end{itemize}

\end{proof}

Using these isomorphisms, one can give another description of the
moduli stack of hyperelliptic curves.
\begin{cor}(\cite[cor. 4.7]{AV})\label{explicit}
The stack $\mathcal{H}_g$ of hyperelliptic curves of genus $g$ can be
represented by:
\begin{itemize}
\item[(i)] If $g$ is even
$$\mathcal{H}_g=[\Asm/GL_2]$$
with action given by $A\cdot f(x)={\rm det}(A)^gf(A^{-1}x)$.
\item[(ii)] If $g$ is odd
$$\mathcal{H}_g=[\Asm/(\mathbb{G}_m\times PGL_2)]$$
with action given by $(\alpha,[A])\cdot f(x)=\alpha^{-2}{\rm det}(A)^{g+1}f(A^{-
1}x)$.
\end{itemize}
\end{cor}
Using this description of $\mathcal{H}_g$ as a quotient, Arsie and Vistoli were
able
to compute
the Picard group of it. For later reference, we include here their instructive
proof. 

First, recall the notion of a functorial Picard group of a stack, as defined by 
Mumford \cite{Mum}, see also \cite[pag. 624]{EG}. For an algebraic stack 
$\mathcal{F}$, 
an element $\E\in {\rm Pic}(\mathcal{F})$ consists of two sets of data:  
\begin{itemize}
\item[(i)] An invertible sheaf $E(\pi)\in \Pic(S)$ for every morphism $S\to 
\mathcal{F}$;
\item[(ii)] For each diagram
$$\xymatrix{
S_1\ar[dr]_{\pi_1} \ar[rr]^f& & S_2\ar[dl]^{\pi_2}\\
&\mathcal{F}&\\
}$$ 
a functorial isomorphism
$\E(f):\E(\pi_1)\stackrel{\cong}{\to}f^*(\E(\pi_2))$.
\end{itemize}
The product of two elements $\E$ and $\E'$ is the line bundle that associates
to every morphism $\pi:S\to \mathcal{F}$ the line bundle $\E(\pi)\otimes 
\E'(\pi)$. 

For a quotient stack $[X/G]$, there is an isomorphism between the 
equivariant Picard group ${\rm Pic}^G(X)$ and the functorial Picard 
group ${\rm Pic}([X/G])$ (see \cite[prop. 18]{EG}). Explicitly, to 
$E\in {\rm Pic}^G(X)$ we associate the element $\E\in {\rm Pic}([X/G])$
whose value on a morphism $S\to [X/G]$, that is a principal $G$-bundle $B\to S$ 
together with an equivariant map $B\to X$ (this is the definition of the 
quotient stack, 
see \cite[section 5.1]{EG}), is just the image of $E$ under the maps ${\rm 
Pic}^G(X)\to
{\rm Pic}^G(B)\cong {\rm Pic}(S)$. Equivalently, viewing an element of 
${\rm Pic}^G(X)$ as a vector bundle $E\to X$ of rank $1$ together with a 
compatible action of $G$, the element $\E\in {\rm Pic}([X/G])$ is the functor 
represented by the line bundle stack $\mathbb{E}:=[E/G]$ over the stack $[X/G]$, 
i.e. we have the following picture:
$$\xymatrix{
E\ar@(ul,ur)^G \ar[r]\ar[d]&X\ar@(ul,ur)^G\ar[d]\\
\mathbb{E}:=[E/G]\ar[r]&[X/G].\\
}$$

\begin{teo}(Arsie-Vistoli, \cite[theo. 5.1]{AV})\label{picstack}
Assume that ${\rm char}(k)$ doesn't divide $2g+2$.
The Picard group of $\mathcal{H}_g$ is
$$
\Pic(\mathcal{H}_g)=
\begin{sis}
\Z/(4g+2)\Z   & \mbox{ if } g \mbox{ is even}\\
\Z/2(4g+2)\Z & \mbox{ if } g \mbox { is odd }.\\
\end{sis}$$
\end{teo}
\begin{demo}
Recall that the equivariant Picard group of a $k$-linear representation
$V$ of a group $G$ is equal to (see \cite[lemma 2]{EG}):
\begin{equation}\label{affine}
\Pic^G(V)=\Pic^G({\rm Spec}(k))=G^*,
\end{equation}
where $G^*={\rm Hom}(G,\G)$ (${\rm Hom}$ is taken in the category of algebraic 
groups). 
Indeed, the group $G$ acts trivially on the automorphism group 
$H^0(V,\O^*_V)=k^*$ of the trivial line bundle on $V$, thus the 
$G$-linearizations of the trivial line bundle on $V$ are elements of ${\rm 
Hom}(G,\G)$. 

Consider the action of $\mathbb{G}_m$ on $\A$, given by
$\alpha\cdot f(x)=\alpha^{-2}f(x)$, and the usual exact sequence:
$$\Z\langle\Delta\rangle\to \Pic^{\mathbb{G}_m}(\A)\to
\Pic^{\mathbb{G}_m}(\Asm)\to 0,$$ where $\Delta=\A-\Asm$ is the
locus defined by the vanishing of the discriminant. More precisely,
here and below $\Delta$ denotes the generator of the subgroup of
such linearizations of the trivial line bundle over $\A$, that
become isomorphic to the trivial linearization over
$\A\backslash\Delta$. Explicitly, since $\Delta$ is an irreducible
hypersurface of degree $4g+2$ (see lemma \ref{Delta}), its equation
is a polynomial $F$ of degree $4g+2$. Multiplication by $F$ defines
an automorphism of the trivial line bundle over
$\A\backslash\Delta$. Since $F(\alpha\cdot f)=\alpha^{-2(4g+2)}F(f)$
for $\alpha\in\G$ and $f\in \A$, we see that the latter automophism
sends a trivial linearization to the trivialization, which
corresponds to the character $-2(4g+2)$ (see also \cite{AV}), and
therefore:
\begin{equation}\label{G_m}
\Pic^{\G}(\Asm)=\Z/2(4g+2)\Z.
\end{equation}

When $g$ is even, consider the two compatible actions:
$$\xymatrix{
\A\ar@(dl,dr)\ar[r]^{id}& \A\ar@(dl,dr)\\
\G\ar@{^{(}->}[r] &GL_2,\\
}$$
where the first action is as before and the second one is that of corollary
\ref{explicit}: $A\cdot f(x)={\rm det}(A)^g f(A^{-1}x)$. Since $GL_2^*\cong \Z$
generated by the determinant morphism, the diagonal inclusion
$\G\hookrightarrow GL_2$ induces a map $GL_2^*\to \G^*$ which is the
multiplication by $2$, i.e. it holds:
\begin{equation}\label{even}
\Pic^{GL_2}(\A)\stackrel{\cdot 2}{\longrightarrow} \Pic^{\G}(\A).
\end{equation}
From the two usual exact sequences
$$\xymatrix{
\Z\langle\Delta\rangle\ar[r]\ar[d]^{id}& \Pic^{GL_2}(\A)\ar[r]\ar[d]^{\cdot
2}&\Pic^{GL_2}(\Asm)\ar[d]\ar[r]&0\\
\Z\langle\Delta\rangle\ar[r] &\Pic^{\mathbb{G}_m}(\A)\ar[r]&
\Pic^{\mathbb{G}_m}(\Asm)\ar[r]& 0,\\
}$$
combined with formula \ref{G_m}, one deduces:
$$\Pic(\mathcal{H}_g)=\Pic^{GL_2}(\Asm)=\Z/(4g+2)\Z,$$
as desired in the even case.

When $g$ is odd, consider the two compatible actions:
$$\xymatrix{
\A\ar@(dl,dr)\ar[r]^{id}& \A\ar@(dl,dr)\\
\G\ar@{^{(}->}[r] &\G\times PGL_2,\\
}$$
where the first action is as before while the second is (according to corollary
\ref{explicit}): $(\alpha,[A])\cdot f(x)=\alpha^{-2}{\rm det}(A)^{g+1}f(A^{-
1}x)$.
In this case, since $(\G\times PGL_2)^*\cong \G^*$, we have an isomorphism
\begin{equation}\label{odd}
\Pic^{\G\times PGL_2}(\A)\stackrel{\cong}{\longrightarrow}\Pic^{\G}(\A).
\end{equation}
Hence from the exact sequences
$$\xymatrix{
\Z\langle\Delta\rangle\ar[r]\ar[d]^{id}& \Pic^{\G\times
PGL_2}(\A)\ar[r]\ar[d]^{\cong}&
\Pic^{\G\times PGL_2}(\Asm)\ar[d]\ar[r]&0\\
\Z\langle\Delta\rangle\ar[r] &\Pic^{\mathbb{G}_m}(\A)\ar[r]&
\Pic^{\mathbb{G}_m}(\Asm)\ar[r]& 0,\\
}$$
combined with formula \ref{G_m}, one deduces:
$$\Pic(\mathcal{H}_g)=\Pic^{\G\times PGL_2}(\Asm)=\Z/2(4g+2)\Z,$$
as desired in the odd case.
\end{demo}
\begin{rem}
This result still holds if one consider the stack $\mathcal{H}_{1,1}$
of families of elliptic curves with a section. In that case the Picard group was
computed by Mumford in his legendary paper \cite{Mum}, and it is isomorphic
to $\Z/12\Z$ (see also \cite[remark 5.5]{AV}).
\end{rem}

Observe that the stack $\mathcal{H}_g^0$ of hyperelliptic curves of genus 
$g$ without extra-automorphisms is isomorphic to the quotient
$$\mathcal{H}_g^0=[\Asm^0/(GL_2/\mu_{g+1})],$$
being equal to the fiber product $\mathcal{H}_g\times_{H_g}H_g^0$. Here $\Asm^0$ 
is the set of forms such that the corresponding $2g+2$ points on $\P^1$ have no 
automorphisms. 

\begin{cor}\label{picstack2}
Assume that ${\rm char}(k)$ doesn't divide $2g+2$.
The Picard group of $\mathcal{H}_g^0$ is
$$
\Pic(\mathcal{H}_g^0)=
\begin{sis}
&\Z/(4g+2)\Z   & \mbox{ if } g \mbox{ is even and } g\neq 2\\
&\Z/5\Z & \mbox{ if } g=2\\
&\Z/2(4g+2)\Z & \mbox{ if } g \mbox { is odd }.\\
\end{sis}$$
\end{cor}
\begin{proof}
For $g>2$ it holds that $\Pic(\mathcal{H}_g)=\Pic(\mathcal{H}_g^0)$ 
since the difference between $\Asm$ and $\Asm^0$ is of codimension at least 2 by 
proposition \ref{dimension}. On the other hand, 
$\Pic(\mathcal{H}^0_2)$ is the quotient of  $\Pic(\mathcal{H}_2)=\Z/10\Z$ 
over the subgroup generated by the divisorial component of $\Asm-\Asm^0$ which, 
in view of lemma \ref{g=2} and the explicit calculations from the proof of 
theorem \ref{picstack}, is the subgroup generated by the residue of $5$
in $\Z/10\Z$. 
\end{proof}

Now we are going to give an explicit description of the generators of 
${\rm Pic}(\H)$ using the functorial description of the Picard group.

\begin{teo}\label{gen}
A generator of ${\rm Pic}(\H)$ is the element $\GG$ that 
associates to a family of hyperelliptic
curves $\pi:\F\to S$ with Weierstrass divisor $W$ the line bundle on $S$:
$$\GG(\pi)=\begin{sis}
\pi_*\left(\omega_{\F/S}^{g+1}\left(-(g-1)W\right)\right)& \text{ if } g \text { 
is even, }\\
\pi_*\left(\omega_{\F/S}^{\frac{g+1}{2}}\left(-\frac{g-1}{2}W\right)\right) & 
\text { if } g \text { is odd. }\\
\end{sis}$$  
\end{teo}
\begin{proof}
From the proof of theorem $\ref{picstack}$, it follows that ${\rm 
Pic}^{GL_2/\mu_{g+1}}
(\Asm)$ is a cyclic group generated by the trivial line bundle $\Asm\times k$ on 
which 
$GL_2/\mu_{g+1}$ acts via a generator of its group of characters. Let us choose 
as the generator the character ${\rm det}^{-(g+1)}$ if $g$ is even and 
${\rm det}^{-\frac{g+1}{2}}$ 
if $g$ is odd (see lemma \ref{quotientgroup}). 
Note that this is true without any assumption on ${\rm char}(k)$ (apart from the 
usual ${\rm char}(k)\neq 2$), while 
the hypothesis that ${\rm char}(k)$ doesn't divide $2g+2$ 
is necessary to compute the order of the Picard group.

We have to express this generator as an element of ${\rm Pic}(\H)$ from the 
point of view of Mumford's functorial description.
Consider the following diagram of $(GL_2/\mu_{g+1})$-equivariant maps 
(the notation is that of theorem \ref{stack}):
$$\xymatrix{
\widetilde{\H}\times k \ar[r]^(0.35){\cong}\ar[d] & \Asm\times k\ar[d]\\
\widetilde{\H}\ar[r]^(0,35){\cong} & \Asm.
}$$

The functor $\widetilde{\H}\times k$ associates to a $k$-scheme $S$ the set
$$\left(\widetilde{\H}\times k\right)(S)=\left\{\C\stackrel{p}{\to} S, \L, 
\L^{\otimes 2}
\stackrel{i}{\hookrightarrow} \O_{\C}
, \phi:(\C,\L)\cong (\P^1_S,\O_{\P^1_S}(-g-1)),\M\right\},$$
where $\M=\O_S$ is the structure sheaf, on which the action of 
$(GL_2/\mu_{g+1})(S)$ is defined via multiplication by ${\rm det}^{-(g+1)}$ if 
$g$ is even 
and ${\rm det}^{-\frac{g+1}{2}}$ if $g$ is odd. 

Let $\P^1_S=\P(V_S)$, where $V$ is a two-dimensional vector space over the 
ground field $k$. From the Euler exact sequence for the 
trivial family $p_S:\P^1_S\to S$ 
$$
0\to \O_{\P^1_S}\to p_S^*(V_S^*) (1)\to 
\omega_{\P^1_S/ S}^{-1}\to 0
$$ 
one deduces a $(GL_2/\mu_{g+1})(S)$-equivariant isomorphism
\begin{equation}\label{Euler}
p_S^*(({\rm det}V_S)^{-1})\otimes \O_{\P^1_S}(2)\cong
\omega^{-1}_{\P^1_S/S},
\end{equation}
where we consider the canonical actions of $(GL_2/\mu_{g+1})(S)$ on $\P^1_S$ and 
on the invertible sheaves involved.
Using projection formula, the fact that $(p_S)_*(\O_{\P^1_S})=\O_S$ and the 
$(GL_2/\mu_{g+1})(S)$-equivariant identity $(\det V_S)^{g+1}=\M$ we get 
$(GL_2/\mu_{g+1})(S)$-equivariant isomorphisms 
$$\begin{sis}
\M&\cong 
(p_S)_*\left(\omega_{\P^1_S/S}^{g+1}\otimes \O_{\P^1_S}(2g+2)\right) 
&\text { if } g \text { is even, }\\
\M&\cong 
(p_S)_*\left(\omega_{\P^1_S/S}^{\frac{g+1}{2}}\otimes \O_{\P^1_S}(g+1)\right) 
&\text{ \: if } g \text{ is odd. }
\end{sis}$$
Now remark that $\phi:(\C,\L)\cong (\P^1_S,\O_{\P^1_S}(-g-1))$ induces a 
canonical isomorphism $\omega_{\C/S}\cong\omega_{\P^1_S/S}$ by the 
$\phi_0$-component. 
Hence the line bundle quotient 
$\mathbb{G}=\left[\left(\widetilde{\H}\times 
k\right)/\left(GL_2/\mu_{g+1}\right)\right]$ 
over $\H'$ is isomorphic to
$$
\mathbb{G}(S)=
\begin{sis}
&\left\{\C\stackrel{p}{\to} S, \L, \L^{\otimes 2}
\stackrel{i}{\hookrightarrow} \O_{\C}, p_*\left(\omega_{\C/S}^{g+1}\otimes 
\L^{-2}\right)\right\} &\text { if } g \text { is even, }\\
&\left\{\C\stackrel{p}{\to} S, \L, \L^{\otimes 2}
\stackrel{i}{\hookrightarrow} \O_{\C}, 
p_*\left(\omega_{\C/S}^{\frac{g+1}{2}}\otimes 
\L^{-1}\right)\right\} &\text { if } g \text { is odd.}\\
\end{sis}
$$
To express the preceding line bundles as push-forward of line bundles on the 
hyperelliptic family $\pi:\F\to S$,  we use formulas \ref{f2} and \ref{f4} 
together with the fact that the line bundles $\omega_{\C/S}^{g+1}\otimes 
\L^{\otimes(-2)}$ and $\omega_{\C/S}^{\frac{g+1}{2}}\otimes 
\L^{-1}$ for $g$ odd are trivial on each fiber of $p$, and we get
$$
\begin{sis}
&f^*\left(\omega_{\C/S}^{g+1}\otimes \L^{-2}\right)=
\omega_{\F/S}^{g+1}\left(-(g-1)W\right) &\text { if } g \text { is even, }\\
&f^*\left(\omega_{\C/S}^{\frac{g+1}{2}}\otimes \L^{-1}\right)=
\omega_{\F/S}^{\frac{g+1}{2}}\left(-\frac{g-1}{2}W\right)& \text{ if } g \text{ 
is odd. }\\
\end{sis}
$$ 
Hence the line bundle $\mathbb{G}$ over $\H$ is equal to
$$\mathbb{G}(S)=
\begin{sis}
&\left\{\F\to S , \pi_*\left(\omega_{\F/S}^{g+1}\left(-(g-
1)W\right)\right)\right\}
& \text{ if } g \text{ is even, }\\
&\left\{\F\to S , \pi_*\left(\omega_{\F/S}^{\frac{g+1}{2}}\left(-\frac{g-
1}{2}W\right)\right)\right\}& 
\text{ if } g \text{ is odd, }\\
\end{sis}
$$
from which the conclusion follows.
\end{proof} 
We can now look at other natural elements of $\Pic(\H)$ and express them in term
of the generator found above. Recall that 
given a family $\pi:\F\to S$ of hyperelliptic curves, there are two natural line 
bundles 
over $\F$: the relative canonical sheaf $\omega_{\F/S}$ and the line bundle 
associated to the Weierstrass divisor $W=W_{\F/S}$. Hence we can consider a 
linear 
combination of them $\omega_{\F/S}^a\otimes \O_{\F/S}(bW)$ and note that 
it restricts on every fiber $F$ of the family to 
$$
\omega_{\F/S}^a\otimes \O_{\F/S}(bW)|_{F}=aK_F+bW_F=a(g-1)g^1_2+b(g+1)g^1_2=
[(a+b)g+(b-a)]g^1_2,
$$
where we used that on a hyperelliptic curve $F$ the canonical class $K_F$
is $(g-1)$-times the unique $g^1_2$ while the Weierstrass divisor $W_F$ is 
$(g+1)$-times the $g^1_2$. Let's call $m(a,b):=[(a+b)g+(b-a)]$ and let's 
consider only those integers $a$ and $b$ for which $m(a,b)\geq 0$. 

Moreover, since on a hyperelliptic curve $F$ it holds that $h^0
(F,\O_F(k g^1_2))=k+1$ for $k\ge 0$,
the push-forward $\pi_*(\omega_{\F/S}^a\otimes \O_{\F/S}(bW))$ is a vector
bundle of rank $m(a,b)+1$ on the base $S$ (see \cite[cor. 12.9]{Har}). Hence we 
can 
define an element $\T$ of ${\rm Pic}(\H)$ by 
$$\T(\pi)={\rm det }\left(\pi_*(\omega_{\F/S}^a\otimes \O_{\F/S}(bW))\right)\in 
{\rm Pic}(S).$$

\begin{teo}\label{elements}
In terms of the generator $\GG$ of $\,{\rm Pic}(\H)$ (see theorem \ref{gen}), 
if $0\leq m(a,b)< g+1$ the element $\T$ is equal to
$$\T=\begin{sis}
\GG^{\frac{(a+b)(m(a,b)+1)}{2}}& \text{ if } g \text{ is even, }\\
\GG^{(a+b)(m(a,b)+1)} & \text{ if } g \text{ is odd, }
\end{sis}$$
and if $m(a,b)\ge g+1$ the element $\T$ is equal to
$$\T=\begin{sis}
\GG^{\frac{(a+b-1)(m(a,b)-g)}{2}}& \text{ if } g \text{ is even, }\\
\GG^{(a+b-1)(m(a,b)-g)} & \text{ if } g \text{ is odd. }
\end{sis}$$

\end{teo}
\begin{proof}[I Proof]
The proof consists in pulling-back the element $\T$ to ${\rm 
Pic}(\widetilde{\H})$ 
and then compare it with the pull-back of $\GG$ as $(GL_2/\mu_{g+1})$-linearized 
sheaves (the notation are the same as in the proofs of theorems
\ref{stack} and \ref{gen}).
 
First of all, we want to express $\T$ as an element $\T'$ of the Picard group of 
$\H'\cong \H$.
Since, by formulas \ref{f2} and \ref{f4}, it holds
$$
f^*\left(\omega_{\C/S}^a\otimes \L^{-(a+b)}\right)
=\left(\omega_{\F/S}^a\otimes \O_{\F/S}(bW)\right), 
$$ 
the element $\T'$ in ${\rm Pic}(\H')$ will associate to  
$\{p:\C\to S, \L, i:\L^{\otimes 2}\hookrightarrow \O_{\C}\}\in \H'(S)$ the 
element
\begin{equation}\label{pushforward}
\T'(S)={\rm det} \:p_*\left(\omega_{\C/S}^a\otimes \L^{-(a+b)}\right)\otimes
{\rm det} \:p_*\left(\omega_{\C/S}^a\otimes \L^{-(a+b)+1}\right)\in 
{\rm Pic}(S).
\end{equation}
Here we used that $f_*(\O_{\mathcal{F}})=\O_{\mathcal{C}}\oplus\L$.
In order to compute the pull-back $\widetilde{\T}$ of $\T'$ to ${\rm 
Pic}(\widetilde{\H})$
we use the isomorphism $\phi:(\C,\L)\cong (\P^1_S,\O(-g-1))$ and the Euler 
formula 
\ref{Euler}, which give
$$
\omega_{\C/S}^a\otimes \L^{-(a+b)}\cong p_S^*((\det V_S)^a)\otimes \O_{\P^1_S}
(-2a)\otimes 
\O_{\P^1_S}\left((a+b)(g+1)\right)=
$$
$$
=p_S^*((\det V_S)^a)\otimes \O_{\P^1_S}(m(a,b)),
$$
and, analogously, 
$$
\omega_{\C/S}^a\otimes \L^{-(a+b)+1}\cong 
p_S^*((\det V_S)^a)\otimes \O_{\P^1_S}(m(a,b)-(g+1)),
$$
where $\P^1_S=\P(V_S)$.\\
Now we take the push-forward through the map $p_S$ and take the determinant, 
obtaining
$$
\det\: (p_S)_*\left(p_S^*(\det V_S)^a)\otimes \O_{\P^1_S}(m(a,b))\right)=
\det\left((\det V_S)^{a} \otimes {\rm Sym}^{m(a,b)}(V_S)\right)=
$$
$$
={(\det V_S)}^{(a+b)(g+1)\frac{(m(a,b)+1)}{2}},$$
where we used the relation $\det({\rm Sym}^n(V_S))=(\det 
V_S)^{\frac{n(n+1)}{2}}$. As for the second sheaf, the push forward is zero if 
$m(a,b)< g+1$. Otherwise the analogous computation leads to the following
$$
\det\: (p_S)_*\left(p_S^*(\det V_S)^a)\otimes \O_{\P^1_S}(m(a,b)-(g+1))\right)=
(\det V_S)^{(a+b-1)(g+1)\frac{m(a,b)-g}{2}}.
$$

Now we conclude recalling from theorem \ref{gen} that the pull-back 
$\widetilde{\GG}$
of the generator $\GG$ to the Picard group of $\widetilde{\H}$ is 
$(\det V_S)^{g+1}$
for $g$ even and $(\det V_S)^{\frac{g+1}{2}}$ for $g$ odd.

\end{proof}
\begin{proof}[II Proof]
We will also show another way to find the expression in terms of the canonical 
generators, which is more explicit and doesn't involve the stack
description of theorem \ref{stack}. 

We use the same notations as above. In addition, let $\tau$ denote the 
invertible ``generator'' sheaf
$\pi_*\left(\omega_{\mathcal{F}/S}^{g+1}(-(g-1)W)\right)$ on the base $S$, and 
let $\varepsilon$ denote the invertible ``generator'' sheaf 
$\pi_*\left(\omega_{\mathcal{F}/S}^{\frac{g+1}{2}}(-\frac{g-1}{2}W)\right)$ for 
the case $g$ odd, so that $\varepsilon^2=\tau$. 

The idea is to express the sheaves $\omega_{\mathcal{C}/S}$ and $\L$ in terms of 
$p^*\tau$ (or $p^*\varepsilon$ for $g$ odd) and a certain invertible sheaf 
$\mathcal{E}$ on $\mathcal{C}$, whose determinant of the direct image via $p$ 
can be expressed in terms of $\tau$ (or $\varepsilon$ for $g$ odd). Then one 
conludes by projection formula, using the relation \ref{pushforward}. 

Suppose $g$ is odd. We claim that in this case $\L\cong
p^*(\varepsilon^{-1})\otimes\omega_{\mathcal{C}/S}^{\frac{g+1}{2}} $. Indeed, 
by \ref{pushforward} $f_*(\omega_{\mathcal{F}/S}^{\frac{g+1}{2}})(-\frac{g-
1}{2}W)=
\omega_{\mathcal{C}/S}^{\frac{g+1}{2}}\otimes(\O_{\mathcal{C}}\oplus\L^{-1}).
$ Thus $\varepsilon=p_*(\omega_{\mathcal{C}/S}^{\frac{g+1}{2}}\otimes\L^{-
1})$, and we get the desired statement, since 
$\omega_{\mathcal{C}/S}^{\frac{g+1}{2}}\otimes\L^{-1}$ is isomorphic to the 
structure sheaf on each fiber of $p$.
Moreover, in lemma \ref{trivpushforward} after the proof of theorem
\ref{D-pic}, we will prove that 
$\det\:p_*(\omega_{\mathcal{C}/S}^{m})$ is a trivial line bundle on the base 
$S$ for any $m\in\Z$. 
Thus we are done in the case $g$ odd, taking 
$\mathcal{E}=\omega_{\mathcal{C}/S}$.

Now we treat the case $g$ even. In this case $\mathcal{C}$ is in fact the 
projectivization of a two-dimensional vector bundle $p_*(\M)$ (see theorem 
\ref{G_1^2}, (ii)), where $\M=\omega_{\mathcal{C}/S}^{\frac{g}{2}}\otimes
\L^{-1}$. Hence from Euler exact sequence on each fiber of the family 
$p:\mathcal{\C}\to S$ we get an exact sequence
$$
0\to\O_{\mathcal{C}}\to p^*(p_*(\M)^*)\otimes\M\to
\omega_{\mathcal{C}/S}^{-1}\to 0.
$$
Thus we get an isomorphism
$$
\det (p^*(p_*(\M)^*))\cong \omega_{\mathcal{C}/S}^{-1}\otimes\M^{-2},
$$
so
$$
\det(p_*(\M))\cong\tau.
$$
Moreover, there are expressions 
$$
\omega_{\mathcal{C}/S}\cong 
p^*\tau\otimes\M^{-2},
$$
$$
\L\cong(p^*\tau)^{\frac{g}{2}}\otimes\M^{-(g+1)}.
$$
Therefore we can take $\mathcal{E}=\M$ for $g$ even.
\end{proof}

Among the elements $\T$ one is of particular interest, namely the Hodge 
line bundle that in our notation is $\mathcal{T}_{1,0}(\F\to S)
=\det \: \pi_*(\omega_{\F/S})$. It is known that, over the complex numbers,
the Hodge line bundle generate the Picard group of 
$\mathcal{M}_g$ (see \cite{AC}). For hyperelliptic curves we have the 
following
\begin{cor}\label{Hodge}
In terms of the generator $\GG$ of $\Pic(\H)$, the Hodge line bundle is equal
to
$$ det \: \pi_*(\omega_{\F/S})=
\begin{sis}
&\GG^{g/2}& \text{ if } g \text{ is even, }\\
&\GG^{g}&\text{ if } g \text{ is odd. }
\end{sis}$$
In particular it generates $\Pic(\H)$ if $g$ is not divisible by $4$
while otherwise it generates a subgroup of index $2$.
\end{cor}   
For $g=2$, this was proved by Vistoli in \cite{Vis} (he computed
the Chow ring of $\mathcal{H}_2=\mathcal{M}_2$ proving that it's 
generated by the Chern classes of the Hodge bundle).

Note also that for $g$ even there is another interesting generator of the 
Picard group of $\H$, that is $\mathcal{T}_{\frac{g}{2},1-\frac{g}{2}}$
(for which it holds that $a+b=1$ and $m(a,b)=1$). The interest
of it is that it is the determinant of the push-forward of a globally 
defined $g^1_2$ (which is very far from being unique!) on the family $\F\to S$, 
that in fact, as we know from section 3, exists in general only for $g$ even.

\section{Comparison between stack and coarse moduli space of hyperelliptic 
curves}

Now we want to compare the stack $\H$ with its coarse moduli scheme
$H_g$ (as well as the open 
substack $\H^0$ with the open subvariety $H_g^0$) .

We introduce a new moduli functor that is "intermediate" between $\H$ and $H_g$. 

\begin{defi}\label{newfunctor}
The moduli functor $\mathcal{D}_{2g+2}$ is the contravariant
functor
$$\mathcal{D}_{2g+2}:Sch_{/k}\rightarrow Set$$
which associates to every $k$-scheme $S$ the set
$$\mathcal{D}_{2g+2}(S)=\left\{
\begin{aligned}
&\C\rightarrow S \text{ family of } \mathbb{P}^1 \text{ and }
D\subset \C
\text{ an effective Cartier divisor }\\
&\text{ finite and \'etale over } S \text{ of degree }
2g+2
\end{aligned}
\right\}_{/\cong}$$ and $\mathcal{D}_{2g+2}^0$ is the subfunctor
of families of effective divisors on $\mathbb{P}^1$ without automorphisms.
\end{defi}
Being without automorphisms means for an effective divisor on $\mathbb{P}^1$
that there are no projective transformation of $\P^1$ that preserves
the divisor. By the results of Lonstead-Kleiman 
(see theorem \ref{hyperfamilies} and \ref{Weierstrass}) 
it follows that there is
a natural transformation of functors $\Psi:\mathcal{H}_g\to
\mathcal{D}_{2g+2}$. Moreover, since over an algebraically closed
field a hyperelliptic curve is uniquely determined (up to
isomorphism) by the $2g+2$ points on $\P^1$ (up to isomorphism)
over which the double cover of $\P^1$ is ramified, it follows that both
these moduli functors have $H_g$ as a coarse moduli scheme. We end
up with the following diagram:
$$\xymatrix{
\mathcal{H}_g \ar[rr]^{\Psi} \ar[rd]_{\Phi_{\mathcal{H}}} & &
\mathcal{D}_{2g+2}
\ar[ld]^{\Phi_{\mathcal{D}}}\\
& {\rm Hom}(-,H_g).& }$$ 

Now we want to prove that $\D$ is an algebraic stack providing a description 
of it as a quotient stack.

\begin{teo}\label{D-stack}
$\D$ is an algebraic stack isomorphic to the quotient
stack \\
$[\Bsm/PGL_2]$, where the action is given by 
$[A]\cdot [f(x)]=[f(A^{-1}\cdot x)]$. Moreover it holds the following 
isomorphism of stacks $[\Bsm/PGL_2]\cong [\Asm/(GL_2/\mu_{2g+2})]$,
with the same action as before. 
\end{teo}
\begin{demo}
We prove the first part of the theorem with a 
strategy analogous to that of theorem \ref{stack} of
Arsie and Vistoli, namely first rigidifying the functor so that it becomes 
a scheme and then viewing this rigidified functor as a principal bundle over 
the original one for the action of a suitable group.

Here the rigidified functor is the functor $\widetilde{\D}$ that associates 
to a $k$-scheme $S$ the set
$$\widetilde{\D}(S)=\{\C\to S, D, \phi:\C\cong \P^1_S\}$$
where $\C\to S$ is a family of $\P^1$, $D$ is an effective Cartier divisor 
as the one in definition \ref{newfunctor} and $\phi$ is an isomorphism
between the family $\C\to S$ and the trivial family $\P^1_S=S\times \P^1_S$.

This rigidified functor is isomorphic to $\Bsm$ (thought as the functor
${\rm Hom}(-,\Bsm)$). In fact an effective smooth divisor of degree $2g+2$  
on $\P^1$ is an element of $\Bsm$ and hence a divisor $D$ 
on $\C\cong \P^1_S$ as above can be identified with an element of $\Bsm(S)$.  

The group sheaf $\underline{\rm Aut}(\P^1)\cong \underline{PGL_2}$ acts on
$\widetilde{\D}$ by composing with the isomorphism $\phi$ and 
it's easy to see that the corresponding action of $PGL_2$ on $\Bsm$ 
is the one given in the statement.

Finally, descent theory implies that the forgetfull morphism 
$\widetilde{\D}\to \D$ makes $\widetilde{\D}$ into a 
$\underline{\rm Aut}(\P^1)$-principal bundle over $\D$, from which
one gets the description of $\D$ as a quotient stack $[\Bsm/PGL_2]$.

To prove the second part of the theorem, observe that, applying lemma  
\ref{quotientgroup}(ii) with $g+1$ replaced by $2g+2$, one deduce an 
isomorphism $GL_2/\mu_{2g+2}\cong \G\times PGL_2$. Moreover
one can check that, under this isomorphism, the corresponding action of 
$GL_2/\mu_{2g+2}\cong \G\times PGL_2$ on $\Asm$ is given by
$(\alpha, [A])\cdot f(x)= \alpha^{-1}\cdot ({\rm det}A)^{g+1}f(A^{-1}\cdot x)$. 
Hence the stack quotient of $\Asm$ by $GL_2/\mu_{2g+2}\cong \G\times PGL_2$
can be taken in two steps: first take the quotient over the subgroup 
$\G/\mu_{2g+2}\cong\G$, 
which is isomorphic to $\Bsm$ since the action is free, and then take the 
quotient over 
$GL_2/\G\cong PGL_2$ with the usual action.
\end{demo}

In the next theorem we compute the Picard group of
$\D$ and compare it with the Picard group of $\H$.
We use the first description of $\D$ in the preceding theorem, although 
everything can be proved also using the second description 
in a spirit similar to theorems \ref{picstack} and  \ref{gen}. 

\begin{teo}\label{D-pic}
Assume that ${\rm char}(k)$ doesn't divide $2g+2$. Then 
$\Pic(\D)=\Z/(4g+2)\Z$ generated by the element $\Gdiv$ that associates
to a family $p:\C\to S$ of $\P^1$ togheter with a Cartier divisor $D$
(as in the definition \ref{newfunctor}) the line bundle on $S$
$$\Gdiv(p)=p_*\left(\omega_{\C/S}^{g+1}(D)\right).$$
Moreover, the natural map $\Pic(\D)\to \Pic(\H)$ is injective and hence 
it's an isomorphism if $g$ is even while it's an inclusion of index $2$ 
if $g$ is odd.  
\end{teo} 
\begin{demo}
In view of the explicit description of theorem \ref{D-stack}, it holds 
that $\Pic(\D)=\Pic^{PGL_2}(\Bsm)$ and we already proved 
(see formula \ref{4.2})
that this a cyclic group of order $4g+2$ generated by
$\O_{\Bsm}(1)$ with its natural $PGL_2$-linearization.
  
To prove the functorial description of $\Gdiv$, we can pull-back this element to 
$\Pic(\widetilde{\D})$ (see theorem \ref{D-stack}) and hence we reduce 
to show the isomorphism of the 
corresponding $PGL_2$-equivariant invertible sheaves for the case when $S=\Bsm$, 
$\mathcal{C}=\Bsm\times\P^1$, and $D$ is the incidence divisor. As in the proof 
of theorem \ref{gen}, from Euler formula \ref{Euler} applied to the given family 
one deduces a $PGL_2$-equivariant isomorphism 
$$
p^*((\det E)^{-(g+1)})\otimes\O_{\P(E)}(2g+2)\cong\omega^{-g-1},
$$
where $\omega$ denotes the (trivial) relative canonical sheaf for the morphism 
$p:\P^1\times\Bsm\to\Bsm$ and $E=V\times\Bsm$ is a trivial two-dimensional 
vector 
bundle on $\Bsm$ such that $\C=\P(E)$. So after taking push-forwards we get
$$
{\rm Sym}^{2g+2}(E)\cong p_*(\omega^{-(g+1)})\otimes(\det E)^{g+1}.
$$
Remark that the group $PGL_2$ acts trivially on $\det E$, thus we get an exact 
sequence of 
$PGL_2$-equivariant sheaves on $\P^1\times\Bsm$ 
$$
0\to p^*\O_{\B}(-1)\to\omega^{-(g+1)}\to\omega^{-(g+1)}|_D\to 0.
$$
So on the base $\Bsm$ there is an equality of 
$PGL_2$-equivariant sheaves $\O_{\B}(-1)\cong p_*(\omega^{-(g+1)}(-D))$ (we use 
the fact that the restriction of $\omega^{-(g+1)}(-D)$ on each fiber is 
trivial), and we get the desired statement.

Finally to study the map $\Pic(\D)\to \Pic(\H)$, 
let us first remark that
$$\Pic^{GL_2/\mu_{g+1}}(\A)=\Pic^{GL_2/\mu_{g+1}}(\AO),$$
since the origin in $\A$ is of codimension $\ge 2$
(see \cite[sect. 2.4, lemma 2]{EG}).
Now consider the compatible diagram:
\begin{equation}\label{basicdiagram}
\xymatrix{
\AO\ar@(dl,dr)_{\G}\ar[r]& \B\ar@(dl,dr)_{\G},\\
}\end{equation}
where the action on the left is given by
$\alpha\cdot f(x)=\alpha^{-2}f(x)$ while on the right is the
trivial one. There is an isomorphism $\Pic^{\G}(\B)\cong\Z\oplus\Z$, where the
first component is generated by $\O_{\B}(1)$ with the trivial action of $\G$,
and the second component is just $\G^*$. The trivialization of the pull-back of
$\O_{\B}(1)$ on $\AO$ is given by the section $f(x)\mapsto(f(x),f(x))\in
\AO\times\A$, so the pull back of the trivial $\G$-linearization of $\O_{\B}(1)$
corresponds to the character $-2\in\Z=\G^*$, because of the action
$\alpha\cdot f(x)=\alpha^{-2}f(x)$. Thus we see, that the composition
$\Z=\Pic^{PGL_2}(\B)\to\Pic^{\G}(\B)\to\Pic^{\G}(\AO)=\Z$ is equal to
multiplication by $-2$.

We can complete the diagram (\ref{basicdiagram}) from above as follows:
$$\xymatrix{
\AO\ar@(ul,ur)^{GL_2/\mu_{g+1}}\ar[r]\ar[r]& \B\ar@(ul,ur)^{PGL_2}\\
\AO\ar@(dl,dr)_{\G}\ar[r]\ar[u]_{id}& \B\ar@(dl,dr)_{\G}\ar[u]_{id}.\\
}$$

So, in the case $g$ even we see from (\ref{even}) that the morphism
$\Pic^{PGL_2}(\B)\to\Pic^{GL_2/\mu_{g+1}}(\AO)=\Pic^{GL_2}(\AO)$ is an
isomorphism $\Z\stackrel{-1}\longrightarrow\Z$, while for $g$ odd we see from
(\ref{odd}) that $\Pic^{PGL_2}(\B)\to\Pic^{GL_2/\mu_{g+1}}(\AO)=\Pic^{\G\times
PGL_2}(\AO)$ is equal to multiplication by $-2$.

Now the conlcusion follows since $\Pic(\D)$ is the quotient of
$\Pic^{PGL_2}(\B)=\Z$ of order $4g+2$ while $\Pic(\mathcal{H}_g)$ 
is the quotient of
$\Pic^{GL_2}(\AO)=\Z$ of order $4g+2$ if $g$ is even and order 
$2(4g+2)$ if $g$ is odd (see theorem \ref{picstack}).
\end{demo}

With the same technique of above, we can prove the following lemma that was used 
in the second proof of theorem \ref{elements}.
\begin{lem}\label{trivpushforward}
For a $\P^1$-family $p:\mathcal{C}\to S$ with an effective Cartier divisor 
$D\subset\mathcal{C}$ etal\'e and finite of degree $2g+2$ over $S$, 
the line bundle $\det \left(p_*\omega_{\mathcal{C}/S}^m\right)$ is 
trivial for any $m\in \Z$.
\end{lem}
\begin{proof}
By pulling-back to $\widetilde{\D}$ (see theorem \ref{D-stack}), 
one reduce to consider the $PGL_2$-equivariant line bundle 
$\det(p_*\omega^m)$ for the trivial family $p:\P(E)\to \Bsm$ together 
with the incidence divisor $D$, where $E=V\times\Bsm$ is a trivial two-
dimensional 
vector bundle. Using Euler formula \ref{Euler}, one expresses 
$\det(p_*\omega^m)$ as a power of $\det(E)$ and hence $PGL_2$ 
acts trivially on it.
\end{proof}

Using theorem \ref{D-pic}, it is possible to 
proof a weaker form of theorem \ref{gen} without computations for stacks. 
Namely, it is possible to proof the statement of theorem \ref{gen} for $g$ even 
and only up to 2-torsion (which is isomorphic to $\Z/2\Z$) for $g$ odd. 

In notations of theorem \ref{gen} the generator $\mathcal{G}$ of the Picard 
group $\Pic(\H)$ corresponds to the residue class of $-1$ in 
Arsie-Vistoli description as a cyclic group (see theorem \ref{picstack}). 
Moreover, from 
the proof of theorem \ref{D-pic} it follows that the map of cyclic groups 
$\Pic(\D)\to\Pic(\H)$ is multiplication by $-1$ for $g$ even and 
by $-2$ for $g$ odd. Thus we have just to reinterpret the generator sheaf 
$\Gdiv$ from 
theorem  \ref{D-pic} in terms of the family of hyperelliptic curves 
and to take the square root in the case $g$ odd. As in the proof of theorem 
\ref{gen} one obtains that
$$
f^*(\omega_{\mathcal{C/S}}^{g+1}(D))
\cong\omega_{\mathcal{F}/S}^{g+1}(-(g+1)W)\otimes\O_{\mathcal{F}}(2W)
=\omega_{\mathcal{F}/S}^{g+1}(-(g-1)W),
$$
where $f:\mathcal{F}\to\mathcal{C}$ is the quotient over the hyperelliptic 
involution, $W\subset \mathcal{F}$ is the Weierstrass divisor and 
$D\subset\mathcal{C}$ is the branch divisor. 


Now we want to study how much 
the stacks $\H$ and $\D$ are far to be finely represented by 
their coarse moduli scheme $H_g$. Since the existence of automorphisms 
is always one of the most seriuos obstruction to the finess of 
moduli scheme, it's very natural to restrict to the open subset $H_g^0$
of hyperelliptic curves without extra-automorphisms as well as to the 
corresponding stacks $\H^0$ and $\D^0$. In fact we get a positive answer
for the stack $\D^0$.

\begin{teo}\label{universaldivisor}
$H_g^0$ is a fine moduli scheme for the functor
$\mathcal{D}_{2g+2}^0$, i.e. the natural transformation of
functors
$$\Phi_{\mathcal{D}}:\mathcal{D}_{2g+2}^0\xrightarrow{\cong}
{\rm Hom}(-,H_g^0)$$
is an isomorphism.
\end{teo}
\begin{proof}
We have to construct a family of $\P^1$ over $H_g^0$ plus an effective
Cartier divisor (finite and \'etale of degree $2g+2$ over
$H_g^0$) that is universal for the functor $\D^0$. To do
this, we consider over $(\S-\Delta)^0$ (the open subset of
$(2g+2)$-tuples without automorphisms) a trivial family of $\P^1$
together with the tautological divisor above it:
$$\xymatrix{
(\S-\Delta)^0\times \P^1 \ar[d] & \mathbb{D}_{2g+2}=\{(D,x): x\in D\}
\ar@{_{(}->}[l]\\
(\S-\Delta)^0. & \\
}$$ Now $PGL_2$ acts (naturally) on $(\S-\Delta)^0$ and diagonally
on  $(\S-\Delta)^0\times \P^1$ and this action clearly preserves
the tautological divisor $\mathbb{D}_{2g+2}$. Moreover, since we
restrict over $(\S-\Delta)^0$, the action of $PGL_2$ is free on
the family of $\P^1$ as well on the divisor $\mathbb{D}_{2g+2}$.
Hence everything passes to the quotient giving the required
universal family of $\P^1$ plus the divisor:
$$\xymatrix{
((\S-\Delta)^0\times \P^1)/PGL_2=\C_g \ar[d] &
D_{2g+2}=\mathbb{D}_{2g+2}/PGL_2
 \ar@{_{(}->}[l]\\
(\S-\Delta)^0/PGL_2=H_g^0. & \\
}$$
Now since the parameterized objects are really without automorphisms, the
existence of a tautological family shows that in fact it's a universal one.
\end{proof}
In view of this result, we can reinterpret the last assertion in 
theorem \ref{D-pic} as follows:
\begin{cor}\label{picmap}
Assume that ${\rm char}(k)$ doesn't divide $2g+2$.
The natural map  $\Pic(H_g^0)\to \Pic(\mathcal{H}_g^0)$ is injective.
Hence it's an isomorphism for $g$ even, 
while it's an inclusion of index $2$ for $g$ odd.
\end{cor}

\begin{rem}
Compare these results with the analogous ones for the moduli spaces of curves 
of genus $g\geq 3$ (results that up to now are known only over the complex 
numbers).
In that case there is an inclusion
$$\Pic(M_g)\hookrightarrow \Pic(\mathcal{M}_g)\cong\Pic(\mathcal{M}_g^0)
\cong\Pic(M_g^0)\cong\Cl(M_g^0)\cong\Cl(M_g).$$ 
It is known that $\Pic(\mathcal{M}_g)\cong\Z$ generated by the Hodge class
(see \cite{Har} and \cite{AC}) but it's still unknown the index of the first 
group into 
the second (see \cite[section 4]{AC}). 
\end{rem}

\begin{rem}
Using theorem \ref{universaldivisor} we could prove lemma \ref{trivpushforward} 
for families of divisors on $\P^1$ of degree $2g+2$ without automorphisms and 
also find the generator of $\Pic(H_g^0)$ repeating the proof of theorem 
\ref{D-pic} saying nothing about stacks. Thus without stack theory we could 
prove a weaker form of theorem \ref{gen} as in the discussion after lemma 
\ref{trivpushforward}, and also theorem \ref{elements} for families of 
hyperelliptic curves without 
extra-automorphisms. In fact, this makes difference only for $g=2$ by 
proposition \ref{dimension}.
\end{rem}

Now we can study the natural trasformation 
$\Phi_{\mathcal{H}}:\H^0\to {\rm Hom}(-, H_g^0)$, 
or in other words we study
how many families of hyperelliptic curves 
(without extra-automorphisms) can have
the same modular map. We generalize Mumford's arguments from the case
of elliptic curves to the case of hyperelliptic curves without 
extra-automorphisms
(see \cite[pag. 49-53, pag. 60-61]{Mum}).
\begin{teo}\label{doublecovers}
Given a map $\phi:S\to H_g^0$, the set of families of
hyperelliptic curves having $\phi$ as a modular map, if non empty,
is a principal homogeneous space for $H^1_{\acute et}(S,\Z/2\Z)$.
\end{teo}
\begin{proof}
Fix a map $\phi:S\to H_g^0$ and
suppose that it is a modular map for some family of hyperelliptic
curves over $S$. Let's denote by $\mathcal{H}_g(S)_{\phi}$ the
(non-empty) set of families of hyperelliptic curves
over $S$ having $\phi$ as modular map. \\
We are going to define an action of $H^1_{\acute et}(S,\Z/2\Z)$ on
$\mathcal{H}_g(S)_{\phi}$ as follows: for a family $\pi:\F\to S$
in $\mathcal{H}_g(S)_{\phi}$ and an element of $H^1_{\acute
et}(S,\Z/2\Z)$ (i.e. a double \'etale cover $f: S'\to S$), we
define a new family $f\cdot \pi: \F'\to S$ of
$\mathcal{H}_g(S)_{\phi}$ by mean of the following diagram:
\begin{equation}\label{diagram}
\xymatrix{
& S'\times_{S} \mathcal{F}  \ar[dl] \ar[dr] \ar[dd] \ar@(ul,ur)[]^{j\times i}&
\\
\mathcal{F}'= S' \times_{S} \mathcal{F}/_{(j\times i)}\ar[dd]_{f\cdot \pi} &
\ar @{} [dr] |{\Box} &\mathcal{F} \ar[dd]^{\pi} \ar@(u,r)[]^i\\
&S' \ar[dr]^f \ar[dl]  \ar@(dl,dr)[]_j& \\
S=S'/j & & S\\
}\end{equation}
where $j$ is the involution on $S'$ that exchanges the two
sheets of the covering $f$, $i$ is the global hyperelliptic
involution (see \ref{hyperfamilies}(i)) and $j\times i$ is the
involution on the fiber product. So the new family $\F'\to S$ is
obtained first by doing the pull-back of the family $\F\to S$ to $S'$
and then by taking the quotient with respect to $i\times j$. Note that also the
first part of the diagram is cartesian and that the original family
$\F\to S$ can be re-obtained by taking the quotient of $S'\times_S \F'$
with respect to the involution $j\times {\rm id}$.

By construction,
over a geometric point $s\in S$ the fibers of $\pi:\F\to S$ and
$f\cdot \pi:\F'\to S$ are the same so that the new family is an
element of $\mathcal{H}_g(S)_{\phi}$ and the definition is
well-posed.

We have to show that this action is simply transitive, namely that given
two families $\pi_1:\F_1\to S$ and $\pi_2:\F_2\to S$ of hyperelliptic curves
in $\mathcal{H}_g(S)_{\phi}$ there exists a unique \'etale double cover of
$S$ which realizes the construction in diagram \ref{diagram}. By general results
of Grothendieck (see \cite{Gro}), there exists a scheme
${\rm Isom}({\pi_1,\pi_2})$
over $S$ whose fiber over the geometric point $s\in S$ is
$${\rm Isom}({\pi_1,\pi_2})_s={\rm Isom}({\pi_1^{-1}(s),\pi_2^{-1}(s)})$$
and hence, since the fibers of our families are hyperelliptic curves
without extra-automorphisms, this is a double \'etale cover of $S$.
Moreover the two families become isomorphic above
${\rm Isom}(\pi_1,\pi_2)$ and the corresponding diagram:
$$\xymatrix{
& \F  \ar[dl] \ar[dr] \ar[dd] & \\
\F_1 \ar[dd]_{\pi_1} \ar @{} [dr] |{\Box}& \ar @{} [dr] |{\Box} &\F_2
\ar[dd]^{\pi_2} \\
&{\rm Isom}(\pi_1,\pi_2)\ar[dr] \ar[dl]  & \\
S & & S\\
}$$
satisfies exactly the property of diagram \ref{diagram}
(see \cite[pag. 61]{Mum}).
Moreover this is the unique double cover with that property
(see \cite[pag. 61]{Mum}).
\end{proof}

\begin{rem}
This proof is in fact an explicit form of a general principle and could be
formulated shorter in the following way. Any two hyperelliptic families with the
same madular map to $H_g^0$ are locally isomorphic in the \'etale topology since
in a suitable \'etale neighborhood both families have a section. Moreover, the
automorphism group of hyperelliptic families is $\Z/2\Z$, so the set of all
families over the base $S$ with the same modular map form a principal homogenous
space over $H^1_{\acute e t}(X,\Z/2\Z)$.
\end{rem}

The non-uniqueness of a family with a given modular map may be also seen
explicitly from the construction of double covers. By theorem
\ref{universaldivisor}, $\mathcal{H}_g(S)_{\phi}$ is the set of families
of hyperelliptic curves that are double covers of $p:\phi^*(\C_g)\to S$ branched
along
$\phi^*(D_{2g+2})$. By the general theory of cyclic covers
(see \cite{Par} or \cite{AV}), such double covers are determined by a line
bundle
$\L$ over $\phi^*(\C_g)$ and an isomorphism
$\L^{\otimes 2}\cong \O_{\phi^*(\C_g)}(-\phi^*(D_{2g+2}))$.

Clearly we may change $\L$ by an $\L'$ such that
$\L\otimes \L'^{-1}=p^*(\mathcal{M})$
for some $\mathcal{M}\in \Pic(S)_2$. We also may change the isomorphism,
multiplying it by a representative of a class from
$\O(\phi^*(\C_g))^*/(\O(\phi^*(\C_g))^*)^2=
\O(S)^*/(\O(S)^*)^2$. The relation to what was said before is provided by the
exact sequence
$$
1\to \O(X)^*/(\O(X)^*)^2\to H^1_{\acute e t}(X,\Z/2\Z) \to \Pic(X)_2\to 0,
$$
which follows from the Kummer exact sequence of sheaves in the \'etale topology.

There is another, stack theoretical, interpretation of theorem 
\ref{doublecovers} and also of the fact that $\Phi_{\mathcal{H}}$ is not 
surjective and its relation to the Brauer group. Namely, we use the second
description of $\mathcal{D}_{2g+2}$ from theorem \ref{D-stack}.

Let us recall that if a group scheme $G$ acts on the scheme $X$ then there is an 
``exact sequence'' of fibered categories
$$
X\to [X/G]\to H^1_{\acute e t}(-,G),
$$
where the action of $G$ carries the groupoid structure on (the category 
associated to) $X$. Instead of groupoid $[X/G](S)$ consider the set of 
equivalency classes of its objects. We will denote it by the same letter.

For any scheme $S$ there is an ``exact sequence'' of sets
$$
X(S)/G(S)\to [X/G](S)\to H^1_{\acute e t}(S,G),
$$
where exactness in the middle term is with respect to the pointed set structure 
on $H^1_{\acute e t}(S,G)$. Explicitly, any element in $[X/G](S)$ is given by an 
\'etale covering $S=\cup_{\alpha}U_{\alpha}$ and maps $f_{\alpha}:U_{\alpha}\to 
X$, $g_{\alpha\beta}:U_{\alpha}\cap U_{\beta}\to G$ such that on the 
intersection $U_{\alpha}\cap U_{\beta}$ there is 
$f_{\alpha}=g_{\alpha\beta}f_{\beta}$ for all $\alpha,\beta$. Moreover, an 
equivalence (i.e. an isomorphisms of the initial groupiod) between 
$\{U_{\alpha},f_{\alpha},g_{\alpha\beta}\}$ and 
$\{U'_{\alpha},f'_{\alpha},g'_{\alpha\beta}\}$ is given by a common subcovering 
$V_{\gamma}$ and maps $g_{\gamma}:V_{\gamma}\to G$ such that after suitable 
restriction $f_{\gamma}=g_{\gamma}f'_{\gamma}$ and $g_{\gamma\delta}=
g^{-1}_{\gamma}g'_{\gamma\delta}g_{\delta}$ in evident notations. Thus the last 
map is just projection to $g_{\alpha\beta}$. 

Now suppose we are given a central extension of group schemes
$$
1\to K\to G\to H\to 1
$$ 
and a free action of $H$ on the scheme $X$. Then for each scheme $S$ the group 
$H^1_{\acute e t}(S,K)$ acts naturally on $[X/G](S)$ by formula 
$\{U_{\alpha},f_{\alpha},g_{\alpha\beta}\}\mapsto
\{U_{\alpha},f_{\alpha},g_{\alpha\beta}k_{\alpha\beta}\}$, where 
$k_{\alpha\beta}$ is a 1-cocycle. 

\begin{pro}\label{stacks}
The set $(X/H)(S)$ is a quotient of $[X/G](S)$ under the action of $H^1_{\acute 
e t}(S,K)$. 
\end{pro}
\begin{proof}
Suppose $\{U_{\alpha},f_{\alpha},g_{\alpha\beta}\}$ and 
$\{U'_{\alpha},f'_{\alpha},g'_{\alpha\beta}\}$ are equivalent in $(X/H)(S)$. By 
definition there exists a smaller subcovering $V_{\gamma}$ and 
$h_{\gamma}:V_{\gamma}\to H$ such that $f_{\gamma}=h_{\gamma}f'_{\gamma}$. 
Taking, if necessary, a smaller subcovering, we may suppose that for each 
$\gamma$ there is $g_{\gamma}:V_{\gamma}\to G$ such that it naturally maps to 
$h_{\gamma}$. Thus multiplying by $g_{\gamma}$ we see that we could suppose from 
the very beginning that $f_{\alpha}=f'_{\alpha}$. Hence we obtain 
$f_{\alpha}=g_{\alpha\beta}f_{\beta}=g'_{\alpha\beta}f_{\beta}$. This means that 
$g_{\alpha\beta}=k_{\alpha\beta}g'_{\alpha\beta}$ for a certain 1-cocycle 
$k_{\alpha\beta}:U_{\alpha\beta}\to K$ since the action of $H$ on $X$ is free.
\end{proof}
Besides, the map $[X/G](S)\to (X/H)(S)$ is not surjective. The obvious 
cohomological obstruction is provided by the image of a given element from 
$(X/H)(S)$ under the composition $(X/H)(S)\to H^1(S,H)\to H^2(S,K)$.

In our case $G=GL_2/\mu_{g+1}$, $H=GL_2/\mu_{2g+2}$, $K=\Z/2\Z$, $X=\Asm^0$ and 
proposition \ref{stacks} becomes theorem \ref{doublecovers}. The cohomological 
obstruction takes values in $H^2_{\acute e t}(S,\Z/2\Z)$ which is a 
reinterpretation of theorem \ref{G_1^2}. Indeed, using the isomorphism from 
lemma \ref{quotientgroup} we see that for $g$ odd the exact sequence of groups 
in question is
$$
0\to\Z/2\Z\to GL_2\stackrel{(\det(\cdot),[\cdot])}\longrightarrow 
\mathbb{G}_m\times PGL_2\to 1
$$ 
while for $g$ even this is 
$$
0\to\Z/2\Z\to \mathbb{G}_m\times PGL_2\stackrel{(2,1)}\longrightarrow 
\mathbb{G}_m\times PGL_2\to 1,
$$
Thus the exact sequence
$$
0\to\Pic(S)/2\Pic(S)\to H^2_{\acute e t}(S,\Z/2\Z)\to{\mathrm Br}(S)_2\to 0.
$$ 
shows that the triviality of the cohomological obstruction means for $g$ odd 
that a certain divisor should be divisible by two, and for $g$ even, in 
addition, that the $\P^1$-family should be Zariski locally trivial.

\begin{exa}\label{explic}
Let $S={\rm Spec}(k)$. Then $H^1_{\acute e t}(S,\Z/2\Z)=k^*/(k^*)^2$.
If we fix a divisor $D\subset\P^1_k$ over $k$ of degree $2g+2$ then
the set of all hyperelliptic curves over $k$, corresponding to the
pair $(\P^1,D)$, may be described as follows: any such
hyperelliptic curve is locally given by the equation
$$
ay^2=P(x)
$$
where $P(x)$ is some fixed equation of the divisor $D$ on
$\mathbb{A}^1\subset \P^1$ and $a$ corresponds to a class from
$k^*/(k^*)^2$.
\end{exa}

In the last part of this section, we are going to investigate 
the existence of a tautological 
family of hyperelliptic curves over an open subset of $H_g$ (compare
\cite[exercise 2.3]{HM} after having replaced universal with tautological!). 

\begin{teo}\label{tautologicalfamily}
There exists a tautological family of hyperelliptic curves over an open
subset of $H_g$ if and only if $g$ is odd. 
\end{teo}
\begin{proof}
Clearly it's enough to restrict to $H_g^0$.
We proved in theorem \ref{universaldivisor} that over $H_g^0$ there exists a
family $\C_g$ of $\P^1$ plus a divisor $D_{2g+2}$ (finite and
\'etale over $H_g^0$ of degree $2g+2$) that are universal. Hence
if a tautological family of hyperelliptic curves exists over an open subset 
$U\subset H_g^0$, then it has to be a double cover of ${\C_g}_{|U}$ 
branched along $D_{2g+2}$.

Now if $g$ is odd, theorem \ref{G_1^2} (iii) gives the existence of a
tautological family over an open subset of $H_g^0$.  

On the other hand for $g$ even, the non-existence of a tautological family
over any open subset of $H_g^0$ will follow from theorem \ref{G_1^2} (ii) 
once we will prove that the family $\C_g\to H_g^0$ is not Zariski locally 
trivially.
 
Let's consider again the situation of theorem
\ref{4.2}:
$$\xymatrix{
\mathbb{D}_{2g+2} \ar@{^{(}->}[r] &(\S-\Delta)^0\times \P^1 \ar[d]
\ar[r] & \C_g \ar[d]
& D_{2g+2} \ar@{_{(}->}[l]\\
& (\S-\Delta)^0 \ar[r] & H_g^0.  &\\
}$$ First of all from \cite[II ex. 6.1]{Har} one gets:
$$
\Pic((\S-\Delta)^0\times \P^1)=(\Z\cdot
p_1^*(\O(1)))/(4g+2)\Z\oplus \Z\cdot p_2^*(\O(1)).
$$
where $p_1$ and $p_2$ are the projections on the first and on the
second factor. 

To compute the Picard group of $\C_g$, we use again the theory of
\emph{equivariant Picard group}  of Mumford (\cite{GIT}). Note
that in this case the action is free so that actually
$$\Pic(\C_g)=\Pic^{PGL_2}((\S-\Delta)^0\times \P^1).$$
Since the action of $PGL_2$ is diagonal, it holds
$$\Pic^{PGL_2}((\S-\Delta)^0\times
\P^1)=\Pic^{PGL_2}((\S-\Delta)^0)\times
\Pic^{PGL_2}(\P^1).$$ We already proved (see \ref{4.2},
\ref{4.4} and the last part of the proof of theorem \ref{picmoduli}) that:
$$\Pic^{PGL_2}((\S-\Delta)^0)=
\begin{sis}
&\Z/(4g+2)\Z& \text{ if } g\geq 3,\\
&\Z/5\Z& \text{ if } g=2.
\end{sis}
$$
generated by the hyperplane section. As for the action
$\sigma:PGL_2\times \P^1\to \P^1$, we have that
$\sigma^*(\O_{\P^1}(1))=p_1^*(\O_{PGL_2}(1))\otimes
p_2^*(\O_{\P^1}(1))$ and, since $\O_{PGL_2}(1)$ is of $2$-torsion
in $\Pic(PGL_2)$, it follows that only $\O_{\P^1}(2)$ admits a
$PGL_2$-linearization or in other words:
$$\Pic^{PGL_2}(\P^1)=\Z\cdot \O_{\P^1}(2).$$
Therefore
$$
\Pic(\C_g)=
\begin{sis}
&(\Z\cdot p_1^*(\O(1)))/(4g+2)\Z\oplus \Z\cdot p_2^*(\O(2))& \text{ if } g\geq
3\\
&(\Z\cdot p_1^*(\O(1)))/5\Z\oplus \Z\cdot p_2^*(\O(2))& \text{ if } g=2.\\
\end{sis}
$$
Hence, since there doesn't exist a line bundle of vertical degree $1$,
by proposition \ref{Zariskitriviality}(4i) the family $\C_g\to H_g^0$ is 
not Zariski locally trivial.

The non existence of a line bundle of vertical degree 1 on the family $\C_g$
over $H_g^0$ may be also deduced from the universality of this family and from
the existence of any family of divisors without automorphisms of degree $2g+2$
on $\P^1$ which has no divisor of horizontal degree 1. For example, as such
family we could take an open subset of the set of all conics in $\P^2$, on which
the divisor of degree $2g+2$ is defined by the intersection with an irreducible
curve of degree $g+1$ in $\P^2$.
\end{proof}

In the preceeding theorem \ref{tautologicalfamily}, we really need to take 
an open subset of $H_g$ for $g$ odd. Namely, the following is true
\begin{pro}
There doesn't exist a tautological family over all $H_g$ (and neither
over $H_g^0$) for $g$ odd.
\end{pro}
\begin{proof}
Clearly it's enough to restrict to $H_g^0$.
There are two ways to prove this fact. In fact they are rather just two 
different ways of looking at the same situation. 

The first way is to suppose that such family exists and consider the 
corresponding morphism 
$H_g^0\to \mathcal{H}_g^0$. It would imply the existence of a projection for 
the inclusion $\Pic(H_g^0)\hookrightarrow \Pic(\mathcal{H}_g^0)$ but this is 
impossible since the first group is isomorphic to $\Z/(4g+2)$ by theorem 
\ref{picmoduli}
and the second is isomorphic to $\Z/2(4g+2)\Z$ by the 
corollary \ref{picstack2}. Thus we get a contradiction.

The second way is to compute explicitly the class of $D_{2g+2}$ in the Picard 
group $\Pic(\C_g)$ (recall that $D_{2g+2}$ is the universal divisor in the 
family $\C_g\to H_g^0$). Since the fiber of $\mathbb{D}_{2g+2}$
over a point $x\in \P^1$ is a hyperplane in $(\S-\Delta)^0$ and
over a point $D\in (\S-\Delta)^0$ consists of $2g+2$ points of
$\P^1$, the class of $\mathbb{D}_{2g+2}$ in the Picard group 
$\Pic((\S-\Delta)^0\times \P^1)=
(\Z\cdot p_1^*(\O(1)))/(4g+2)\Z\oplus \Z\cdot p_2^*(\O(1))$ is equal to
$$
[\mathbb{D}_{2g+2}]=(\overline{1}, 2g+2).
$$
Thus the first component of the class of $D_{2g+2}$ in the Picard group 
$\Pic(\C_g)$ is still equal to $\overline{1}\in \Z/(4g+2)\Z$ so it is 
undivisible by 2 in the Picard group. Hence by remark \ref{remark} 
there isn't any tautological family over $H_g^0$.
\end{proof}

Note that the situation is different for $g=1$ as the 
following remark shows (in this case $\Pic(H_1^0)=0$). 
\begin{rem}
There exists a tautological family
over $H_1^0\cong\mathbb{A}_j^1-\{0,1728\}$ (the isomorphism is given by
associating
to every elliptic curve its $j$-function). The following is an explicitly
example (it's not unique!) of such a family (see \cite[page 58]{Mum}):
$$y^2=x^3+\frac{27}{4}\cdot\frac{1278-j}{j}(x+1).$$
\end{rem}
\begin{rem}
If one considers the moduli space of "framed" hyperelliptic curves 
(i.e. hyperelliptic curve $C$ plus a fixed
double cover $C\to \P^1$), which is just $\S-\Delta$ without taking the quotient
for $PGL_2$, then one can prove that there doesn't exist a universal 
(neither a tautological!) family above it (see \cite{Ran}).
Nevertheless such a tautological family exists over an open subset: for example,
if we remove the hyperplane consisting of tuples containing the point at 
infinity
then the usual equation $y^2=P(x)$, with $P(x)$ a monic polynomial of degree 
$2g+2$ with distinct roots, defines a tautological hyperelliptic curve
(see \cite{Ran}). 
\end{rem}

\section{Application}

There is an interesting application of the theory developed above. Consider a
family $\F\to S$ of smooth hyperelliptic curves of genus $g$ with automorphism
group $\Z/2\Z$ over a regular irreducible base $S$. Let us make two assumptions:
\begin{itemize}
\item the corresponding modular map $S\to H_g^0$ is dominant and generically
finite,
\item the family $\F$ satisfies the conditions of proposition \ref{criterion}.
\end{itemize}

By its universal property the nontrivial element $\alpha$ in ${\rm
Br}(k(H_g^0))$, which corresponds the restriction of the family $\C_g$ on the
generic point, becomes trivial in ${\rm Br}(k(S))$ (see the cohomological
interpretation of proposition \ref{Zariskitriviality}). Recall the following
well-known fact (see \cite[page 12]{Ser}):
\begin{lem}
If a profinite group $H$ is a subgroup of order $n$ inside a profinite group
$G$, then for any $G$--module $M$ and $i\ge 1$ the composition of restriction
and corestriction maps $H^i(G,M)\stackrel{res}\to H^i(H,M) \stackrel{cor}\to
H^i(G,M)$ is equal to the multiplication by $n$.
\end{lem}

This lemma implies that $[k(S):k(H_g^0)]$ must be even since $\alpha$ is of
order 2 being a class of a conic over $k(H_g^0)$.

\begin{rem}
We cannot give a more precise statement about the divisibility of
$[k(S):k(H_g^0)]$ even if we replace the second assumption by a stronger one:
$\F$ has a rational section. Indeed, we may take any double cover $S$ of $H^0_g$
over which the pull-back of $\C_g$ is Zariski locally trivial and then use the
explicit construction from the last example to define a desired family of
hyperelliptic curves over a Zariski open subset in $S$.
\end{rem}

\begin{rem}
If we replace smooth hyperelliptic curves of genus $g$ with automorphism group
$\Z/2\Z$ by smooth complex curves of genus $g$ without automorphisms and assume
the analogous hypotesis  about a family of such curves (namely that the modular
map is
generically finite and dominant and the family has a rational section), then the
answer
will be that the degree of a
modular map should be always divisible by $2g-2$ (\cite[lemma 5]{cap03}).
It follows from a deep statement that says that, over the complex numbers,
the relative Picard group of the
universal family over $M_g^0$ --- the fine moduli space of smooth curves of
genus $g$ without automorphisms --- is generated by the relative canonical
sheaf (this was first claimed by Franchetta \cite{Fra} but a correct proof is
due
to Harer \cite{Harer} and Arbarello-Cornalba \cite{AC}).
There is an analogous statement about families of trigonal curves (in arbitrary
characteristic)
, which is also proved by considering the relative Picard group of the universal
family
(see \cite{GV}).
However, the proof for hyperelliptic case, as presented here, has a
rather different spirit.
\end{rem}

\end{document}